\documentclass[oneside, a4paper,reqno]{amsart}
\usepackage{pdfsync}
\usepackage{stmaryrd}
\usepackage{mathrsfs}
\usepackage{multicol}
\usepackage{amsmath, amsthm, amscd, amssymb, latexsym, eucal}
\usepackage[all]{xy}
\usepackage{stackrel}
 \calclayout \makeatletter
 \calclayout \makeatletter
\def\serieslogo@{} \def\@setcopyright{} \makeatother
\usepackage{multienum}

\usepackage[colorinlistoftodos]{todonotes}

\usepackage{hyperref}
\usepackage{color}
\usepackage[nospace,noadjust]{cite}

\hypersetup{colorlinks=true,
     breaklinks=true,
     linkcolor=blue,
     citecolor=red,
     bookmarks=true,
     pageanchor=true}

\makeatletter
\renewcommand*\env@matrix[1][c]{\hskip -\arraycolsep
  \let\@ifnextchar\new@ifnextchar
  \array{*\c@MaxMatrixCols #1}}
\makeatother

\usepackage{color}

\numberwithin{equation}{section}
\newtheorem{thm}{Theorem}[section]
\newtheorem*{main-thm}{Main Theorem}
\newtheorem*{thmmain}{Theorem}
\newtheorem{cor}[thm]{Corollary}
\newtheorem{lem}[thm]{Lemma}
\newtheorem{prop}[thm]{Proposition}

\theoremstyle{definition}
\newtheorem{defn}[thm]{Definition}
\newtheorem{rem}[thm]{Remark}
\newtheorem{exam}[thm]{Example}

\newtheorem*{ackn}{Acknowledgment}

\newcommand{\lxr}{\longrightarrow}

\newcommand{\A}{\mathscr A}
\newcommand{\B}{\mathscr B}

\newcommand{\F}{\mathcal F}

\newcommand{\M}{\mathcal M}

\newcommand{\T}{\mathcal T}
\newcommand{\U}{\mathcal U}
\newcommand{\V}{\mathcal V}

\newcommand{\X}{\mathcal X}
\newcommand{\Y}{\mathcal Y}
\newcommand{\Z}{\mathcal Z}

\newcommand{\mt}{\mathsf{T}}

\newcommand{\mU}{\mathsf{U}}
\newcommand{\mz}{\mathsf{Z}}

\newcommand{\mr}{\mathsf{r}}
\newcommand{\mq}{\mathsf{q}}
\newcommand{\mi}{\mathsf{i}}

\newcommand{\ml}{\mathsf{l}}
\newcommand{\me}{\mathsf{e}}

\newcommand{\map}{\mathsf{p}}

\DeclareMathOperator*{\Ker}{\mathsf{Ker}}
 \DeclareMathOperator*{\Image}{\mathsf{Im}}

\DeclareMathOperator*{\Mod}{\mathsf{Mod}-\!}

\DeclareMathOperator*{\End}{\mathsf{End}}
 \DeclareMathOperator*{\smod}{\mathsf{mod}-\!}

\DeclareMathOperator*{\proj}{\mathsf{proj}}

\DeclareMathOperator{\Hom}{\mathsf{Hom}}

\DeclareMathOperator*{\Ext}{\mathsf{Ext}}

  \DeclareMathOperator*{\op}{\mathsf{op}}

   \DeclareMathOperator*{\Ab}{\A\!\textit{b}}

\DeclareMathOperator*{\DMor}{\mathsf{DMor}}

\newcommand{\iso}{\cong}

\newcommand{\iden}{\operatorname{Id}\nolimits}

\newsavebox{\proofbox}
\savebox{\proofbox}{\begin{picture}(7,7)%
  \put(0,0){\framebox(7,7){}}\end{picture}}

\begin{document}

\title[]{Ladders of compactly generated triangulated categories \\ and preprojective algebras}
%Infinite ladders induced by preprojective algebras}

\author[C. Psaroudakis]{Nan Gao and Chrysostomos Psaroudakis
}
{\address{Nan Gao \\ Department of Mathematics, Shanghai University, Shanghai 200444, PR China}
\email{nangao@shu.edu.cn}
\address{Chrysostomos Psaroudakis\\
Institute of Algebra and Number Theory, University of Stuttgart, Pfaffenwaldring 57, 70569 Stuttgart, Germany}
\email{chrysostomos.psaroudakis@mathematik.uni-stuttgart.de}}

\date{\today}

\thanks{Nan Gao is supported by Natural Science foundation of China (11771272).}

\thanks{Chrysostomos Psaroudakis was supported  by the Norwegian Research Council (NFR 221893) under the project \textit{Triangulated categories in Algebra} at the Norwegian University of Science and Technology (NTNU)
}

%\thanks{The authors would like to thank Steffen Koenig, Julian K\"ulshammer, Frederik Marks and Jorge Vit\'oria for useful discussions and valuable comments.}

\keywords{Recollement, Ladder, Preprojective algebra, Compactly generated triangulated category.}

\subjclass[2010]{16E10;16E65;16G;16G50;16S50}

\begin{abstract}
In this paper we characterize when a recollement of compactly generated triangulated categories admits a ladder of some height going either upwards or downwards. As an application, we show that the
derived category of the preprojective algebra of Dynkin type $\mathbb{A}_n$ admits a periodic infinite ladder, where the one outer term in the recollement is the derived category of a differential graded algebra.
\end{abstract}

\maketitle

\section{Introduction}

Ladders of triangulated categories were introduced by Beilinson, Ginzburg and Schechtman \cite{BGS} in their attempt to formalize the equivalence between derived categories of graded modules over symmetric and exterior algebras, known now as Koszul duality. A ladder is a recollement of triangulated categories $(\U,\T,\V)$ together with a (possible infinite) sequence of triangle functors going upwards or downwards such that any three consecutive rows form a recollement of triangulated categories. Recall that a recollement $(\U,\T,\V)$ is both a localization and a colocalization sequence of triangulated categories $\U\lxr \T\lxr \V$. Recollements of triangulated categories were introduced by Beilinson, Bernstein and Deligne in their fundamental work on perverse sheaves \cite{BBD}. Today, recollements play an important role in various contexts since the main idea is to encode information for $\T$ from the possibly simpler outer terms $\U$ and $\V$.

In the context of representation theory of finite dimensional algebras, ladders of derived categories have recently attracted a lot of attention due to their connection with several homological and $K$-theoretic invariants, see \cite{AKLY, ChenXi, HanQin, Happel, Keller:cyclichomology}. Another important aspect of ladders of derived categories of rings is the relation to the problem of lifting or restricting recollements to different levels of derived categories, see \cite[Theorem III]{AKLY}. Moreover, by \cite[Section 5]{AKLY} the height of the ladder measures how far an algebra is from being derived simple \cite{W}, $\rm {i.e.}$ its derived category is not the middle term of a non-trivial recollement with outer terms derived categories of algebras. For example, indecomposable symmetric and selfinjective algebras are derived simple at the level of bounded derived categories, see \cite{LY} and \cite{ChenXi} respectively. A standard example of ladders arises from the derived category of a triangular matrix ring \cite[Example 3.4]{AKLY}. Moreover, in the case that the underlying rings are Gorenstein algebras, it was proved in \cite{ZhangZhangZhouZhu} that the derived category of certain triangular matrix algebras admits an infinite ladder.

In light of the above applications, it is natural to study the existence of ladders of compactly generated triangulated categories and search for more examples. In the case of unbounded derived categories of finite dimensional algebras, Angeleri H\"ugel-K\"onig-Liu-Yang (\cite{AKLY}) studied the existence of ladders. The first aim of this paper is to characterize when a recollement of compactly generated triangulated categories admits a ladder of some height going either upwards or downwards (see Theorem~\ref{thmcompactlygenladder}). This can be regarded as a mild generalization of related results in \cite{AKLY}. Our second aim is to present a new example of a recollement which admits an infinite ladder. More precisely, consider the preprojective algebra $\Pi_n(Q)$ in the sense of Gelfand and Ponomarev \cite{GelfandPonomarev}, where $Q$ is a quiver of Dynkin type $\mathbb{A}_n$. For a finite dimensional algebra $\Lambda$ over a field $k$ denote by $\Pi_n(\Lambda, Q)$ the algebra $\Lambda\otimes_{k}\Pi_n(Q)$. We show the following (see Theorem~\ref{doublemorphisminfiniteladder})$\colon$

\begin{thmmain}
The derived category of $\Pi_n(\Lambda, Q)$ admits an infinite ladder $(\mathsf{D}^{}(\Gamma), \mathsf{D}^{}(\Pi_n(\Lambda, Q)), \mathsf{D}^{}(\Lambda))$ of period four, where $\Gamma$ is a $\mathsf{dg}$ algebra, and this ladder restricts to bounded as well as to perfect complexes.
\end{thmmain}

The key ingredients of the proof are$\colon$$(a)$ a characterization of when a recollement of compactly generated triangulated categories admits a ladder of some height, and $(b)$ a careful analysis of the recollement situation of $\Pi_n(\Lambda, Q)$ at the level of module categories. Interestingly, the periodicity of this infinite ladder is not the consequence of the existence of a Serre functor (in contrast to the examples in \cite[Example 3.6]{AKLY}).
Also, since $\Pi_n(Q)$ is a finite dimensional selfinjective algebra it follows from \cite{CK} that $\Pi_n(Q)$ is derived simple, in contrast to the theorem's infinite ladder at the level of unbounded derived categories.
Note that our infinite ladder involves the derived category of a $\mathsf{dg}$ algebra whereas the definition of derived simplicity only allows derived categories of algebras. Hence, we cannot expect to obtain a better recollement situation for $\Pi_n(\Lambda, Q)$ than the one we have in the above Theorem.

The structure of the paper is as follows. In Section~\ref{sectiontworecoladjoints} we collect preliminary notions and results on recollements of triangulated categories and the existence of adjoint functors. In Section~\ref{sectionladdercompagen} we provide necessary and sufficient conditions for a recollement of compactly generated triangulated categories to admit a ladder of some height going either upwards or downwards. This characterization is proved in Theorem~\ref{thmcompactlygenladder} and Theorem~\ref{thmcompactlygenladderdual}. In Section~\ref{Section4laddersmono} we provide sufficient conditions for a recollement of derived categories of $\mathsf{dg}$ algebras to restrict to a recollement between the full subcategories of $\mathsf{dg}$ modules with finite dimensional total cohomology (Proposition~\ref{propresticrecotodfd}). After recalling the description of the module category of $\Pi_n(\Lambda, Q)$, we show in Proposition~\ref{proprecolpreprojalg} that there is a recollement of module categories where the left term is the module category of $\Pi_{n-1}(\Lambda, Q)$ and the right term is the module category of $\Lambda$. In Theorem~\ref{doublemorphisminfiniteladder} we prove the second main result of this paper as stated above. We show that the ladder of $\Pi_n(\Lambda, Q)$ is related to the Nakayama functor, we also explain why Theorem~\ref{doublemorphisminfiniteladder} works only for Dynkin type $\mathbb{A}_n$, and moreover we discuss the relation of $\Pi_2(\Lambda, Q)$ with certain Morita rings.

\begin{ackn}
The authors would like to thank Steffen Koenig, Julian K\"ulshammer, Frederik Marks and Jorge Vit\'oria for useful discussions and valuable comments. The authors wish to thank the referee for the useful suggestions and remarks.
\end{ackn}

\section{Recollements and adjoint functors}
\label{sectiontworecoladjoints}

We start this section with the notion of a recollement due to Beilinson-Bernstein-Deligne \cite{BBD}.

\begin{defn}
\label{defnrecoltriang}
A {\bf \textsf{recollement}} of triangulated categories, denoted by $(\U,\T,\V)$, is a diagram
\begin{equation}
\xymatrix@C=0.5cm{
\U \ar[rrr]^{\mathsf{i}} &&& \T \ar[rrr]^{\mathsf{e}}  \ar @/_1.5pc/[lll]_{\mathsf{q}}  \ar
 @/^1.5pc/[lll]_{\mathsf{p}} &&& \V
\ar @/_1.5pc/[lll]_{\mathsf{l}} \ar
 @/^1.5pc/[lll]_{\mathsf{r}}
 }
\end{equation}
of triangulated categories and triangle functors satisfying the following conditions$\colon$
\begin{enumerate}
\item[\bf 1.] $(\mq,\mi,\map)$ and $(\ml,\me,\mr)$ are adjoint triples.

\item[\bf 2.] The functors $\mi$, $\ml$, and $\mr$ are fully faithful.

\item[\bf 3.] $\Image{\mi}=\Ker{\me}$, where $\Ker{\me}$ is the full subcategory of $\T$ consisting of objects $X$ with $\me (X)=0$, and $\Image{\mi}$ is the full subcategory of $\T$ consisting of
objects $\mi (X)$ for $X\in \U$.
\end{enumerate}
\end{defn}
In the following, we usually write $\mi(\U)$, $\ml(\V)$ and $\mr(\V)$ for $\Image{\mi}$, $\Image{\ml}$ and $\Image{\mr}$ respectively. Note that Definition~\ref{defnrecoltriang} comes with a pair of triangles arising from the units and counits of the adjoint pairs. In particular, it is easy to see that for any object $X$ in $\T$ we have triangles $\mi\map(X)\lxr X\lxr \mr\me(X)\lxr \mi\map(X)[1]$ and $\ml\me(X)\lxr X\lxr \mi\mq(X)\lxr \ml\me(X)[1]$ in $\T$.

It is well known that recollements of triangulated categories correspond bijectively to torsion, torsion-free triples, \textsf{TTF}-triples for short, in triangulated categories (see \cite[Section 1.4.4]{BBD}, \cite[Remark 2.14]{BR} and \cite[Section 4.2]{Ni}). We explain how this bijection works since it is used later.

Let $(\X,\Y,\Z)$ be a \textsf{TTF}-triple in a triangulated $\T$, i.e. $(\X,\Y)$ and $(\Y,\Z)$ are torsion pairs (\cite[Definition 2.1, Chapter I]{BR}) and $\X$, $\Y$ and $\Z$ are triangulated categories. Note that the latter notion is also known as stable t-structures, see \cite{Miyachi}.
Then we have the adjoint pairs $(\mi_{\X},\mathsf{R}_{\X})$, $(\mathsf{L}_{\Y}, \mi_{\Y})$, $(\mi_{\Y}, \mathsf{R}_{\Y})$ and $(\mathsf{L}_{\Z}, \mi_{\Z})$, as indicated in the following diagrams$\colon$
\[
\xymatrix@C=0.5cm{
\X \ar[rrr]^{\mi_{\X}} &&& \T  \ar @/^1.5pc/[rrr]^{\mathsf{L}_{\Y}}  \ar
 @/^1.5pc/[lll]^{\mathsf{R}_{\X}} &&& \Y \ar[lll]_{\mi_{\Y}} \ \ \ (1)  } \ \ \ \ \ \ \ \ \ \
 \xymatrix@C=0.5cm{
\Y \ar[rrr]^{\mi_{\Y}} &&& \T  \ar @/^1.5pc/[rrr]^{\mathsf{L}_{\Z}}  \ar
 @/^1.5pc/[lll]^{\mathsf{R}_{\Y}} &&& \Z \ar[lll]_{\mi_{\Z}} \ \ \ (2)  }
\]
where $\mi_{\X}$, $\mi_{\Y}$ and $\mi_{\Z}$ are the inclusion functors. Associated with the \textsf{TTF}-triple $(\X,\Y,\Z)$ there is a recollement of triangulated categories as follows$\colon$
\begin{equation}
\label{recolarisingfromTTF}
\xymatrix@C=0.5cm{
\Y \ar[rrr]^{\mi_{\Y}} &&& \T \ar[rrr]^{\mathsf{R}_{\X}} \ar
@/_1.5pc/[lll]_{\mathsf{L}_{\Y}}  \ar
 @/^1.5pc/[lll]_{\mathsf{R}_{\Y}} &&& \X
\ar @/_1.5pc/[lll]_{\mi_{\X}} \ar
 @/^1.5pc/[lll]_{\mi_{\Z}\mathsf{L}_{\Z}\mi_{\X}}
 }
\end{equation}

Clearly $\Image{\mi_{\Y}}=\Ker{\mathsf{R}_{\X}}$ and the functors $\mi_{\X}$, $\mi_{\Y}$ are fully faithful. Thus, it remains to prove that $\mi_{\Z}\mathsf{L}_{\Z}\mi_{\X}$ is the right adjoint of $\mathsf{R}_{\X}$, i.e. there is a natural isomorphism $\Hom_{\X}(\mathsf{R}_{\X}(T),X)\iso \Hom_{\T}(T,\mi_{\Z}\mathsf{L}_{\Z}\mi_{\X}(X))$ for every $T\in \T$ and $X\in \X$. Note that from the adjoint triple $(\mi_{\X},\mathsf{R}_{\X},\mi_{\Z}\mathsf{L}_{\Z}\mi_{\X})$ it follows that the functor $\mi_{\Z}\mathsf{L}_{\Z}\mi_{\X}$ is also fully faithful. First, we have $\Hom_{\X}(\mathsf{R}_{\X}(T),X)\iso \Hom_{\T}(\mi_{\X}\mathsf{R}_{\X}(T),\mi_{\X}(X))$. From the diagram $(1)$ we have the triangle $\mi_{\X}\mathsf{R}_{\X}(T)\lxr T\lxr \mi_{\Y}\mathsf{L}_{\Y}(T)\lxr \mi_{\X}\mathsf{R}_{\X}(T)[1]$. Applying the functor $\Hom_{\T}(-,\mathsf{i}_{\Z}\mathsf{L}_{\Z}\mathsf{i}_{\X}(X))$ and using that $\Hom_{\T}(\mi_{\Y}(\Y),\mi_{\Z}(\Z))=0$ we get the isomorphism
\[
\Hom_{\T}(T,\mi_{\Z}\mathsf{L}_{\Z}\mi_{\X}(X))\iso \Hom_{\T}(\mi_{\X}\mathsf{R}_{\X}(T),\mi_{\Z}\mathsf{L}_{\Z}\mi_{\X}(X)).
\]
Next, from diagram $(2)$ we have the triangle $\mi_{\Y}\mathsf{R}_{\Y}\mi_{\X}(X)\lxr \mi_{\X}(X)\lxr \mi_{\Z}\mathsf{L}_{\Z}\mi_{\X}(X)\lxr \mi_{\Y}\mathsf{R}_{\Y}\mi_{\X}(X)[1]$. Applying the functor $\Hom_{\T}(\mathsf{i}_{\X}\mathsf{R}_{\X}(T),-)$ and using that $\Hom_{\T}(\mi_{\X}(\X),\mi_{\Y}(\Y))=0$ we get the isomorphism
\[
\Hom_{\T}(\mi_{\X}\mathsf{R}_{\X}(T),\mi_{\X}(X))\iso \Hom_{\T}(\mi_{\X}\mathsf{R}_{\X}(T),\mi_{\Z}\mathsf{L}_{\Z}\mi_{\X}(X)).
\]
The above isomorphisms show that $(\mathsf{R}_{\X}, \mi_{\Z}\mathsf{L}_{\Z}\mi_{\X})$ is an adjoint pair. Hence the diagram $(\ref{recolarisingfromTTF})$ is a recollement of triangulated categories.
Conversely, if $(\mathcal{U}, \mathcal{T}, \mathcal{V})$ (diagram~\ref{defnrecoltriang}) is a recollement of triangulated categories it is easy to check that the triple $(\ml(\V),\mi(\U), \mr(\V))$ is a \textsf{TTF}-triple in $\T$.

A sequence $\U\stackrel{\mi}{\lxr} \T\stackrel{\me}{\lxr} \V$ of triangle functors $\mi$ and $\me$ between triangulated categories is said to be
exact if the following four conditions are satisfied: (1) The functor $\mi$ is fully faithful. (2) The composition $\me\mi$ is zero.
(3) ${\rm Im}\mi={\Ker}\me$. (4) The functor $\me$ induces a triangle equivalence between the Verdier quotient of $\T$ by ${\Image}\mi$ and $\V$.
The base of a recollement $(\U,\T,\V)$ is an exact sequence $\U\lxr \T\lxr \V$. The following useful result shows that given an exact sequence as above, it suffices to have adjoints (left and right) either for the functor $\mi\colon \U\lxr \T$ or for the functor $\me\colon \T\lxr \V$ to obtain a recollement situation $(\U,\T,\V)$.

\begin{lem} \textnormal{(\cite[Theorem $1.1$]{CPS}, \cite[Theorem $2.1$]{CPSstratif})}
\label{inducedadjoints}
Let $\U\stackrel{\mi}{\lxr} \T\stackrel{\me}{\lxr} \V$ be an exact sequence of triangle functors. Then the following hold:
\begin{enumerate}
\item The functor $\mi$ admits a left adjoint functor $\mq$ if and only if the functor $\me$ admits a left adjoint functor $\ml$.

\item The functor $\mi$ admits a right adjoint functor $\map$ if and only if the functor $\me$ admits a right adjoint functor $\mr$.
\end{enumerate}
In this case, the functor $\ml$ $($respectively, the functor $\mr)$ is fully faithful.
\end{lem}

As mentioned above, it is useful to know when the inclusion functor of a triangulated subcategory $\U$ in a triangulated category $\T$ admits a left or right adjoint. In the following result we provide necessary and sufficient conditions for the inclusion functor $\mi\colon \U\lxr \T$ to have a right adjoint. The dual statement for left adjoints is left to the reader. We first recall some notions and fix some notations.

A full subcategory $\X$ of an additive category $\A$ is called {\bf contravariantly finite} if for any object $A$ in $\A$ there is a morphism $f\colon X_A\lxr A$ in $\A$ with $X_A$ in $\X$ such that the map
$\Hom_{\A}(X',f)\colon$ $\Hom_{\A}(X',X_A)\lxr \Hom_{\A}(X',A)$ is surjective for every object $X'$ in $\X$. Dually we define a subcategory to be covariantly finite and if it is both covariantly and contravariantly finite then it is called {\bf functorially finite}. Also, if $\X$ is a class of objects in $\A$, then we denote by $\X^{\bot}=\{A\in \A \mid \Hom_{\A}(\X,A)=0\}$ the right orthogonal subcategory of $\X$ and by ${^{\bot}\X}=\{A\in \A \mid \Hom_{\A}(A,\X)=0\}$ the left orthogonal subcategory of $\X$.

\begin{lem}\textnormal{(\cite{KV}, \cite[Chapter I, Proposition 2.3]{BR}, \cite[Lemma 3.1]{Bondal})}
\label{lemcharactforadj}
Let $\T$ be a triangulated category with a triangulated subcategory $\U$. The following statements are equivalent.
\begin{enumerate}
\item The inclusion functor $\mi\colon \U\lxr \T$ has a right adjoint, i.e. there is a functor $\map\colon \T\lxr \U$ such that $(\mi,\map)$ is an adjoint pair.

\item The subcategory $\U$ is contravariantly finite in $\T$ and for every object $X$ in $\T$ there is a triangle $U_X\lxr X\lxr U_X'\lxr U_X[1]$ in $\T$ such that the map $U_X\lxr X$ is a right $\U$-approximation of $X$ in $\T$ and $U'_X$ lies in $\U^{\bot}$.

\item $(\U, \U^{\bot})$ is a stable t-structure in $\T$.
\end{enumerate}
\begin{proof}
(i) $\Longrightarrow$ (ii): Let $X$ be an object in $\T$. Using the adjunction isomorphism of the pair $(\mi,\map)$, it follows that the counit map $\map(X)=\mi\map(X)\lxr X$ is a right $\U$-approximation of $X$ in $\T$. We denote by $U_X$ the object $\map(X)$. Consider now the triangle $(*)\colon U_X\lxr X\lxr U_X'\lxr U_X[1]$ in $\T$ and let $U'$ be an object in $\U$. Applying the functor $\Hom_{\T}(U',-)$ to $(*)$ we get the following long exact sequence$\colon$
\[
\xymatrix@C=0.5cm{
(U',U_X) \ar[r]^{\iso} & (U',X) \ar[r] & (U',U_X')  \ar[r]^{} & (U',U_X[1]) \ar[r]^{\iso  } & (U',X[1])  }
\]
Note that the above isomorphisms follow from the adjunction isomorphism of the adjoint pair $(\mi,\map)$ together with $\mi$ being the inclusion functor. This shows that $\Hom_{\T}(U',U_X')=0$ and therefore $U'_X$ lies in $\U^{\bot}$.

(ii) $\Longrightarrow$ (i): By the assumption there is a right $\U$-approximation $U_X\lxr X$ for every $X$ in $\T$. We claim that the assignment $X\mapsto \map(X):=U_X$ induces a functor $\map\colon \T\lxr \U$ which is a right adjoint of the inclusion functor $\mi$. We first show that the above assignment gives a well defined functor.

Let $g\colon Y\lxr X$ be a morphism in $\T$ and consider a right $\U$-approximation $U_{Y}\lxr Y$ of $Y$ in $\T$ such that $U'_Y$ lies in $\U^{\bot}$.
Since $\Hom_{\T}(U_Y,U'_X)=0$ and $\Hom_{\T}(U_Y,U_X'[-1])=0$, we obtain the following commutative diagram$\colon$
\[
\xymatrix{
U_Y \ar[r]^{} \ar[d]_{h} & Y \ar[d]^{g} \ar[r] & U_Y' \ar[d]  \ar[r]^{} & U_Y[1] \ar[d]  \\
U_X \ar[r]^{} & X \ar[r] & U_X'  \ar[r]^{} & U_X[1]  }
\]
where the morphism $h$ is unique. We now show that the right $\U$-approximation $U_X\lxr X$ is the unique up to isomorphism right $\U$-approximation of $X$ in $\T$. Let $V_X\lxr X$ be another right $\U$-approximation of $X$ in $\T$. Then we have the following commutative diagram and $h'\circ h=\iden_{U_X}\colon$
\[
\xymatrix{
U_X \ar[r]^{h} \ar[d]_{} & V_X \ar[d]^{} \ar[r]^{h'} & U_X \ar[d]   \\
X \ar[r]^{\iden_{X}} & X \ar[r]^{\iden_{X}} & X   }
\]
We infer that the map $h$ is an isomorphism.

Finally, let $U'$ be an object in $\U$. Applying the functor $\Hom_{\T}(U',-)$ to the given triangle $U_X\lxr X\lxr U_X'\lxr U_X[1]$,
we get the isomorphism ${\rm Hom}_{\U}(U',U_X)\cong {\rm Hom}_{\T}(U', X)$. This means that $(\mi,\map)$ is an adjoint pair.

(ii) $\Longleftrightarrow$ (iii): This implication (ii) $\Longrightarrow$ (iii) is clear. Conversely, assume that $(\U, \U^{\bot})$ is a stable t-structure. Then for any object $X\in \T$ there is a triangle
$U_X\lxr X\lxr U_X'\lxr U_X[1]$ with $U_X\in \U$ and $U_X'\in \U^{\bot}$. Let $U$ be an object of $\U$. Applying the functor $\Hom_{\T}(U,-)$ to this triangle, we get that the induced morphism $\Hom_{\T}(U, U_X)\lxr \Hom_{\T}(U, X)$ is an isomorphism. Clearly, the map $U_X\lxr X$ is a right $\U$-approximation.
\end{proof}
\end{lem}

The next result shows that in some cases we can lift adjoint functors from abelian categories to derived categories. Its proof is standard, see for instance \cite{Maltsiniotis}. This is used in the proof of the Main Theorem in Section~\ref{Section4laddersmono}.

\begin{lem}
\label{exactadjoint}
Let $\A$ and $\B$ be abelian categories such that $\mathsf{D}(\A)$ and $\mathsf{D}(\B)$ exists. Assume that there is an adjoint pair of exact functors $(F,G)$, i.e. $\xymatrix{
  F\colon \A  \ar@<-.7ex>[r]_-{} & \B \ {{:G}} \ar@<-.7ex>[l]_-{} }$, between $\A$ and $\B$. Then there is an adjoint pair $(F,G)$ between the unbounded derived categories of $\A$ and $\B$ which restricts also to the bounded derived categories. In particular, if $G\colon \B\lxr \A$ is fully faithful, then the induced functor $G\colon \mathsf{D}(\B)\lxr \mathsf{D}(\A)$ is fully faithful.
\end{lem}

\section{Compactly generated triangulated categories and ladders}
\label{sectionladdercompagen}

Our aim in this section is to characterize when a recollement of compactly generated triangulated categories admits a ladder of some height. We first recall the notion of a ladder due to Beilinson-Ginzburg-Schechtman \cite[Section 1.5]{BGS}, see also \cite[Section $3$]{AKLY}.

\begin{defn}
\label{defnladder}
A {\bf \textsf{ladder}} $(\U,\T,\V)$ is a finite or infinite diagram of triangulated categories and triangle functors$\colon$
\begin{equation}
\label{ladderdefn}
\xymatrix@C=0.5cm{
\U  \ar @/^3.0pc/[rrr]^{\mq^1} \ar @/_3.0pc/[rrr]_{\vdots}^{\map^1} \ar[rrr]^{\mi} &&& \T   \ar @/^3.0pc/[rrr]^{\ml^1} \ar @/_4.5pc/[lll]_{\vdots}  \ar[rrr]^{\me } \ar @/_1.5pc/[lll]_{\mq} \ar @/_3.0pc/[rrr]_{\vdots }^{\mr^1}  \ar @/^1.5pc/[lll]_{\map} &&& \V
\ar @/_1.5pc/[lll]_{\ml} \ar @/_4.5pc/[lll]_{\vdots} \ar
 @/^1.5pc/[lll]_{\mr}   }
\end{equation}
such that any three consecutive rows form a recollement of triangulated categories. Multiple occurrence of the same recollement is allowed. The {\bf \textsf{height}} of the ladder $(\U,\T,\V)$ is the number of recollements contained in it (counted with multiplicities). A ladder is {\bf \textsf{periodic}}, if there exists a positive integer $n$ such that the $n$-th recollement going upwards (respectively, going downwards) in the ladder $(\U,\T,\V)$ is equivalent to the recollement $(\U,\T,\V)$ which is considered to be a ladder of height one. The minimal such positive integer $n$ is the {\bf \textsf{period}} of the ladder.
\end{defn}

Let $\T$ be a triangulated category with small coproducts. An object $X$ in $\T$ is called {\bf \textsf{compact}} if the functor $\Hom_{\T}(X,-)\colon \T\lxr \Ab$ preserves coproducts. The compact objects in $\T$ form a thick triangulated subcategory which we denote by $
\T^{\mathsf{c}}$. Then $\T$ is {\bf \textsf{compactly generated}} if $
\T^{\mathsf{c}}$ is skeletally small and the vanishing $\Hom_{\T}(X,Y)=0$ for all $X\in \T^{\mathsf{c}}$ implies that $Y=0$. In other words, $\T$ is compactly generated if $\T$ admits a set of compact generators. In the following, we say that $(\U, \T, \V)$ is a recollement of compactly generated triangulated categories which we mean that $\U$, $\V$ and $\T$ are compactly generated triangulated categories.

In the sequel we need the following useful results which are consequences of Brown representability for compactly generated triangulated categories.

\begin{lem}\textnormal{(\cite[Theorem $4.1$, Theorem $5.1$]{Neeman})}
\label{adjfromBR}
Let $F\colon\mathcal{S}\lxr \mathcal{T}$ be a triangle functor betweeen compactly generated triangulated categories. Assume that the functor $F$ has a right adjoint $G\colon\mathcal{T}\lxr \mathcal{S}$. Then the following are equivalent$\colon$
\begin{enumerate}
\item The functor $G$ preserves coproducts.

\item There is an adjoint triple $(F,G,H)$.

\item The functor $F$ preserves compact objects.

\end{enumerate}
\end{lem}

Note that in a recollement  $(\U, \T, \V)$ (diagram~\ref{defnrecoltriang}) of compactly generated triangulated categories,  $\ml$ and $\mq$ preserve compact objects.

\begin{lem}\textnormal{(\cite[Lemma 3.2]{Neeman})}
\label{lemdevissage}
Let $\T$ be a compactly generated triangulated category. Let $\mathcal{S}$ be a full coproduct-closed triangulated subcategory containing a set of compact generators of $\T$. Then $\mathcal{S}=\T$.
\end{lem}

Given a recollement $(\U, \T, \V)$ of triangulated categories and triangulated subcategories $\X$, $\Y$ and $\Z$ in $\U$, $\T$ and $\V$ respectively, we say that $(\U, \T, \V)$ restricts to an upper recollement $(\X, \Y, \Z)$ if there is a diagram$\colon$
\[
\xymatrix@C=0.5cm{
\X \ar[rrr]^{\mi} &&& \Y \ar[rrr]^{\mathsf{e}}  \ar @/_1.5pc/[lll]_{\mq}  &&& \Z
\ar @/_1.5pc/[lll]_{\mathsf{l}}  }
\]
such that $\X\stackrel{\mi}{\lxr}\Y\stackrel{\me}{\lxr} \Z$ is an exact sequence of triangle functors where $\mi$ admits a left adjoint $\mq$ and $\me$ admits a left adjoint $\ml$.

We are now ready to characterize when a recollement of compactly generated triangulated categories admits a ladder of height two going downwards. This result generalizes \cite[Proposition~$3.2$, Lemma~$4.3$]{AKLY}.

\begin{prop}
\label{propheighttwogoingdown}
Let $(\U, \T, \V)$ be a recollement of compactly generated triangulated categories. The following statements are equivalent$\colon$
\begin{enumerate}
\item There is a ladder of height two going downwards.

\item $\map$ admits a right adjoint $\map^1$.

\item $\mr$ admits a right adjoint $\mr^1$.

\item $\mi$ preserves compact objects.

\item $\me$ preserves compact objects.

\item $\ml\circ \me$ preserves compact objects.

\item The pair $(\mr(\V), \mr(\V)^{\bot})$ is a stable t-structure in $\T$.

\item The recollement $(\U, \T, \V)$ restricts to an upper recollement $(\U^{\mathsf{c}}, \T^{\mathsf{c}}, \V^{\mathsf{c}})$$\colon$
\[
\xymatrix@C=0.5cm{
\U^{\mathsf{c}} \ar[rrr]^{\mi} &&& \T^{\mathsf{c}} \ar[rrr]^{\mathsf{e}}  \ar@/_1.5pc/[lll]_{\mq}  &&& \V^{\mathsf{c}}
\ar @/_1.5pc/[lll]_{\mathsf{l}}  }
\]
\end{enumerate}
In this case, $\mr(\V)$ is a functorially finite subcategory in $\T$.
\begin{proof} (i) $\Longleftrightarrow$ (ii) $\Longleftrightarrow$ (iii) follows from Definition~\ref{defnladder} and Lemma~\ref{inducedadjoints}.

(ii) $\Longleftrightarrow$ (iv)  and (iii) $\Longleftrightarrow$ (v) follow from Lemma~\ref{adjfromBR}.

(v) $\Longrightarrow$ (vi): This follows from the fact that if the triangle functors $F$ and $G$ between compactly generated triangulated categories preserve compact objects, so does the composition $GF$.

(vi) $\Longrightarrow$ (v): Let $X$ be a compact object in $\T$ and $\{V_i \ | \ i\in I\}$ a set of objects in $\V$. Then the commutativity of the following diagram$\colon$
\[
\xymatrix{
\Hom_{\V}(\me(X),\coprod_{i\in I}V_i) \ar[d]^{\iso}_{\ml \ \text{f.f.}} \ar[rrrr]^{} & & &  &  \coprod_{i\in I}\Hom_{\V}(\me(X),V_i) \ar[d]^{\iso}_{\ml \ \text{f.f.}}   \\
\Hom_{\T}(\ml\me(X),\ml(\coprod_{i\in I}V_i)) \ar[rr]^{\iso}_{\ml \ \text{left adjoint}} && \Hom_{\T}(\ml\me(X),\coprod_{i\in I}\ml(V_i)) \ar[rr]^{\iso}_{\ml\me(X)\in \T^{\mathsf{c}}} && \coprod_{i\in I}\Hom_{\T}(\ml\me(X),\ml(V_i))     }
\]
shows that the object $\me(X)$ is compact in $\V$.

(iii) $\Longleftrightarrow$ (vii): Since the functor $\mr$ has a right adjoint $\mr^1$ if and only if the inclusion functor $\mr(\V)\lxr \T$ has
a right adjoint $\T\lxr \mr(\V)$, the result follows from Lemma~\ref{lemcharactforadj}.

So far we have proved that the first seven conditions are equivalent. We now show (iii) $\Longrightarrow$ (viii). By assumption and the arguments above, $(\mq, \mi)$ is an adjoint pair between $\T^{\mathsf{c}}$ and $\U^{\mathsf{c}}$ where  $\mi\colon \U^{\mathsf{c}}\lxr \T^{\mathsf{c}}$ is fully faithful. Moreover, $(\ml,\me)$ is an adjoint pair between $\V^{\mathsf{c}}$ and $\T^{\mathsf{c}}$ where $\ml\colon \V^{\mathsf{c}}\lxr \T^{\mathsf{c}}$ is fully faithful, and $\mi(\U^{\mathsf{c}})\subseteq \Ker{\me}$. Let $Y$ be an object in $\T^{\mathsf{c}}$ such that $\me(Y)=0$. Then there is an object $U$ in $\U$ such that $\mi(U)=Y$. We claim that $U$ lies in $\U^{\mathsf{c}}$. Let $\{U_i \ | \ i\in J\}$ be a set of objects in $\U$. Then we have the following commutative diagram$\colon$
\[
\xymatrix{
\Hom_{\U}(U,\coprod_{i\in J}U_i) \ar[d]^{\iso}_{\mi \ \text{f.f.}} \ar[rrrr]^{} & & &  &  \coprod_{i\in J}\Hom_{\U}(U,U_i) \ar[d]^{\iso}_{\mi \ \text{f.f.}}   \\
\Hom_{\T}(\mi(U),\mi(\coprod_{i\in J}U_i)) \ar[rr]^{\iso}_{\mi \ \text{left adjoint} \ } && \Hom_{\T}(\mi(U),\coprod_{i\in J}\mi(U_i)) \ar[rr]^{\iso}_{\mi(U)\in \T^{\mathsf{c}}} && \coprod_{i\in J}\Hom_{\T}(\mi(U),\mi(U_i))     }
\]
and therefore the object $U$ is compact in $\U$. Hence $\mi(\U^{\mathsf{c}})=\Ker{\me}$.

Conversely, since the functor $\me$ preserves compact objects by (viii), it follows from Lemma~\ref{adjfromBR} that $\mr$ has a right adjoint $\mr^1$. Thus (iii) holds.

From the dual of Lemma~\ref{lemcharactforadj}, the subcategory $\mr(\V)$ is always covariantly finite in $\T$ since the inclusion functor $\mr(\V)\lxr \T$ has the quotient functor $\me$ as a left adjoint. On the other hand, since the functor $\mr$ has a right adjoint $\mr^1$ if and only if the inclusion functor $\mr(\V)\lxr \T$ has a right adjoint $\T\lxr \mr(\V)$, we get from Lemma~\ref{lemcharactforadj} that $\mr(\V)$ is contravariantly finite in $\T$.  This implies that $\mr(\V)$ is functorially finite in $\T$.
\end{proof}
\end{prop}

We continue with the dual result of Proposition~\ref{propheighttwogoingdown}, that is,  when a recollement of compactly generated triangulated categories admits a ladder of height two going upwards.  Given a recollement $(\U, \T, \V)$, usually one can not expect a new one by interchanging $\U$ and $\V$. Note that there is some difference between Proposition~\ref{propheighttwogoingdown} and Proposition~\ref{propheighttwogoingup}, in particular, compare the upper recollements that we obtain at the level of compact objects.

\begin{prop}
\label{propheighttwogoingup}
Let $(\U, \T, \V)$ be a recollement of compactly generated triangulated categories. The following statements are equivalent$\colon$
\begin{enumerate}
\item There is a ladder of height two going upwards.

\item $\mq$ admits a left adjoint $\mq^1$.

\item $\ml$ admits a left adjoint $\ml^1$.

\item The pair $({^{\bot}\ml(\V)}, \ml(\V))$ is a stable t-structure in $\T$.

\item The recollement $(\U, \T, \V)$ induces an upper recollement $(\V^c, \T^c, \U^c)$$\colon$
\[
\xymatrix@C=0.5cm{
\V^{\mathsf{c}} \ar[rrr]^{\ml} &&& \T^{\mathsf{c}} \ar[rrr]^{\mq}  \ar @/_1.5pc/[lll]_{\ml^1}  &&& \U^{\mathsf{c}}
\ar @/_1.5pc/[lll]_{\mq^1}  }
\]
\end{enumerate}
In this case, $\ml(\V)$ is a functorially finite subcategory in $\T$.
\begin{proof}
(i) $\Longleftrightarrow$ (ii) $\Longleftrightarrow$ (iii) follows from Definition~\ref{defnladder} and Lemma~\ref{inducedadjoints}.

The equivalence of (iii) and (iv) follows from the dual of Lemma~\ref{lemcharactforadj}. Now we show that (iii) $\Longrightarrow$ (v). By assumption and the above arguments, the triangle functors $\ml^1$, $\ml$, $\mq^1$ and $\mq$ preserve compact objects. Hence, we get the diagram
\[
\xymatrix@C=0.5cm{
\V^{\mathsf{c}} \ar[rrr]^{\ml} &&& \T^{\mathsf{c}} \ar[rrr]^{\mq}  \ar @/_1.5pc/[lll]_{\ml^1}  &&& \U^{\mathsf{c}}
\ar @/_1.5pc/[lll]_{\mq^1}  }
\]
such that $(\ml^1, \ml)$ and $(\mq^1, \mq)$ are adjoint pairs. Clearly the functors $\ml$ and $\mq^1$ are fully faithful. It remains to show that $\Image{\ml}=\Ker{\mq}$. This proof follows in the same way as the proof of the corresponding part in Proposition~\ref{propheighttwogoingdown}, the details are left to the reader.

Conversely, assume that (v) holds. We only need to show that the adjunctions $(\ml^1,\ml)$ and $(\mq^1,\mq)$ extend from the categories of compact objects to the whole triangulated categories. More precisely, we show that $(\ml^1,\ml)$ is an adjoint pair between $\V$ and $\T$, i.e. statement (iii) holds. Let $A$ be an object in $\T^{\mathsf{c}}$ and consider the full subcategory of $\V\colon$
\[
{_A\M}=\big\{Y\in \V \ | \ f_{A,Y}\colon \Hom_{\V}(\ml^1(A),Y[k])\stackrel{\iso}{\lxr} \Hom_{\T}(A,\ml(Y[k])), \forall k\in \mathbb{Z} \big\}
\]
It is easy to check that ${_A\M}$ is a triangulated subcategory of $\V$ and by assumption $\V^{\mathsf{c}}$ is contained in ${_A\M}$. Let $(Y_i)_{i\in I}$ be a family of objects in ${_A\M}$. From the commutative diagram
\[
\xymatrix{
\Hom_{\V}(\ml^1(A),\coprod_{i\in I}Y_i) \ar[d]^{\iso}_{\ml^1(A)\in \V^{\mathsf{c}}} \ar[rrrr]^{f_{A,\coprod_{i\in I}Y_i}} & & &  & \Hom_{\T}(A,\ml(\coprod_{i\in I}Y_i)) \ar[d]^{\iso}_{\ml \ \text{left adjoint}}   \\
\coprod_{i\in I}\Hom_{\V}(\ml^1(A),Y_i) \ar[rr]^{\iso}_{Y_i\in {_A\M} \ }  &&  \coprod_{i\in I}\Hom_{\T}(A,\ml(Y_i)) \ar[rr]^{\iso}_{A\in \T^{\mathsf{c}}} &&  \Hom_{\T}(A,\coprod_{i\in I}\ml(Y_i))
}
\]
it follows that the map $f_{A,\coprod_{i\in I}Y_i}$ is an isomorphism and therefore the triangulated subcategory ${_A\M}$ is closed under coproducts.  Then, from Lemma~\ref{lemdevissage} we obtain that $(*)\colon {_A\M}=\V$ for every compact object $A$ in $\T$.  This means that the map $f_{A,Y}$ is an isomorphism for every $A$ in $\T^{\mathsf{c}}$ and $Y$ in $\V$. On the other hand, for an arbitrary but fixed $Y$ in $\V$ consider the following subcategory of $\T\colon$
\[
\M_Y=\big\{X\in \T \ | \ f_{X,Y}\colon \Hom_{\V}(\ml^1(X[k]),Y)\stackrel{\iso}{\lxr} \Hom_{\T}(X[k],\ml(Y)), \forall k\in \mathbb{Z} \big\}
\]
It follows easily as above that $\M_Y$ is a coproduct-closed triangulated subcategory of $\T$ and by the relation $(*)$ we deduce that $\T^{\mathsf{c}}\subseteq \M_Y$. Then, Lemma~\ref{lemdevissage} implies that $\M_Y=\T$ for all $Y$ in $\V$. We infer that $(\ml^1,\ml)$ forms an adjoint pair between the triangulated categories $\V$ and $\T$.

From Lemma~\ref{lemcharactforadj}, the subcategory $\ml(\V)$ is always contravariantly finite in $\T$ since the inclusion functor $\ml(\V)\lxr \T$ has the quotient functor $\me$ as a right adjoint. On the other hand, since the functor $\ml$ has a left adjoint $\ml^1$ if and only if the inclusion functor $\ml(\V)\lxr \T$ has a left adjoint $\T\lxr \ml(\V)$, we get from the dual of Lemma~\ref{lemcharactforadj} that $\ml(\V)$ is covariantly finite in $\T$. We infer that $\ml(\V)$ is functorially finite in $\T$.
\end{proof}
\end{prop}

We are now ready to characterize when a recollement of compactly generated triangulated categories admits a ladder of height $n$ going downwards. Parts of the following result generalizes \cite[Theorem 4.4]{AKLY}.
For simplicity in the presentation of the next result we fix a positive odd integer $n\geq 3$. The case that $n$ is even is treated in a similar way.

\begin{thm}
\label{thmcompactlygenladder}
Let $n\geq 3$ a positive odd integer and $(\U, \T, \V)$ be a recollement of compactly generated triangulated categories
\[
\xymatrix@C=0.5cm{
\U \ar[rrr]^{\mathsf{i}} &&& \T \ar[rrr]^{\mathsf{e}}  \ar @/_1.5pc/[lll]_{\mathsf{q}}  \ar
 @/^1.5pc/[lll]_{\mathsf{p}} &&& \V
\ar @/_1.5pc/[lll]_{\mathsf{l}} \ar
 @/^1.5pc/[lll]_{\mathsf{r}}
 }
\]
The following statements are equivalent:
\begin{enumerate}
\item The recollement $(\U, \T, \V)$ admits a ladder of height $n$ going downwards
\begin{equation}
\label{ladderdefn}
\xymatrix@C=0.5cm{
\U \ar @/_3.0pc/[rrr]_{\vdots}^{\map^1} \ar[rrr]^{\mi} &&& \T    \ar @/^5pc/[lll]^{p^{n-1}}  \ar[rrr]^{\me } \ar @/_1.5pc/[lll]_{\mq} \ar @/_3.0pc/[rrr]_{\vdots }^{\mr^1}  \ar @/^1.5pc/[lll]_{\map} &&& \V
\ar @/_1.5pc/[lll]_{\ml} \ar @/^5pc/[lll]^{r^{n-1}} \ar
 @/^1.5pc/[lll]_{\mr}   }
\end{equation}

\item There are sequences of functors $\map^1, \map^3, \ldots, \map^{n-4}\colon \U\lxr \T$ and $\map^2, \ldots, \map^{n-3}\colon \T\lxr \U$ which preserve compact objects and the pairs $(\map,\map^1), \ldots, (\map^{n-2}, \map^{n-1})$ are adjoint pairs.

\item There are sequence of functors $\mr^1, \mr^3, \ldots, \mr^{n-4}\colon \T\lxr \V$ and $\mr^2, \mr^4, \ldots, \mr^{n-3}\colon \V\lxr \T$ which preserve compact objects and the pairs $(\mr,\mr^1), \ldots, (\mr^{n-2}, \mr^{n-1})$ are adjoint pairs.

\item The recollement $(\U, \T, \V)$ restricts to a recollement $(\U^{\mathsf{c}}, \T^{\mathsf{c}}, \V^{\mathsf{c}})$ which admits a ladder of height $n-2$ going downwards.

\item There is a sequence of triangulated subcategories $\ml(\V)$, $\mi(\U)$, $\mr(\V)$, $\map^1(\U)$, $\mr^2(\V)$, $\map^3(\U), \ldots, \mr^{n-3}(\V)$, $\map^{n-2}(\U)$, $\mr^{n-1}(\V)$ in $\T$ such that
\[
\big(\ml(\V),\mi(\U),\mr(\V)\big), \big(\mi(\U),\mr(\V),\map^1(\U)\big), \ldots, \big(\mr^{n-3}(\V), \map^{n-2}(\U), \mr^{n-1}(\V)\big)
\]
are \textnormal{\textsf{TTF}}-triples in $\T$.

\item  There is a sequence of  \textnormal{\textsf{TTF}}-triples in $\T^{\mathsf{c}}$ of the form$\colon$$(\ml(\V)\cap \T^{\mathsf{c}},\mi(\U)\cap \T^{\mathsf{c}},\mr(\V)\cap \T^{\mathsf{c}})$, $(\mi(\U)\cap \T^{\mathsf{c}},\mr(\V)\cap \T^{\mathsf{c}},\map^1(\U)\cap \T^{\mathsf{c}}), \ldots, (\mr^{n-5}(\V)\cap \T^{\mathsf{c}}, \map^{n-4}(\U)\cap \T^{\mathsf{c}}, \mr^{n-3}(\V)\cap \T^{\mathsf{c}})$.
\end{enumerate}
\begin{proof}
For simplicity we prove the result for $n=3$. The case where $(\U, \T, \V)$ has a ladder of height $n\geq 5$ is treated similarly.

(i) $\Longleftrightarrow$ (ii) $\Longleftrightarrow$ (iii): Assume that $(\U, \T, \V)$ admits a ladder of height three going downwards. Then we have the adjoint triples $(\mi,\map,\map^1)$, $(\map,\map^1,\map^2)$ and $(\me,\mr,\mr^1)$, $(\mr,\mr^1,\mr^2)$. Hence from Lemma~\ref{adjfromBR} we get that the functors $\map$, $\map^1$, $\mr$ and $\mr^1$ preserve compact objects. This implies (ii) and (iii). Conversely,  (ii)$\Longrightarrow$ (i) and (iii)$\Longrightarrow$ (i) follow from Lemma~\ref{adjfromBR}.

(iv) $\Longrightarrow$ (i): Restricting the recollement to compact objects means precisely that the functors $\mi$, $\map$ and $\me$, $\mr$ preserve compact objects. Then Lemma~~\ref{adjfromBR} provides us with the extra adjoints such that the recollement $(\U, \T, \V)$ admits a ladder of height three going downwards.

(i) $\Longrightarrow$ (iv): Since the recollement $(\U, \T, \V)$ admits a ladder of height three going downwards, we have the adjoint triples$\colon$ $(\mi,\map,\map^1)$, $(\map,\map^1,\map^2)$, $(\me,\mr,\mr^1)$ and $(\mr,\mr^1,\mr^2)$. From Lemma~\ref{adjfromBR} the triangle functors $\map$ and $\mr$ preserve compact objects. On the other hand, there is an induced upper recollement by Proposition~\ref{propheighttwogoingdown}
\[
\xymatrix@C=0.5cm{
\U^{\mathsf{c}} \ar[rrr]^{\mi} &&& \T^{\mathsf{c}} \ar[rrr]^{\mathsf{e}}  \ar@/_1.5pc/[lll]_{\mq}  &&& \V^{\mathsf{c}}
\ar @/_1.5pc/[lll]_{\mathsf{l}}  }
\]
It yields a recollement
\[
\xymatrix@C=0.5cm{
\U^{\mathsf{c}} \ar[rrr]^{\mathsf{i}} &&& \T^{\mathsf{c}} \ar[rrr]^{\mathsf{e}}  \ar @/_1.5pc/[lll]_{\mathsf{q}}  \ar
 @/^1.5pc/[lll]_{\mathsf{p}} &&& \V^{\mathsf{c}}
\ar @/_1.5pc/[lll]_{\mathsf{l}} \ar
 @/^1.5pc/[lll]_{\mathsf{r}}
 }
\]
and in this case the height is one.

(i) $\Longrightarrow$ (v): By the assumption we have the following three recollements of triangulated categories$\colon$
\[
\xymatrix@C=0.3cm{
\U \ar[rrr]^{\mathsf{i}} &&& \T \ar[rrr]^{\mathsf{e}}  \ar @/_1.5pc/[lll]_{\mathsf{q}}  \ar
 @/^1.5pc/[lll]_{\mathsf{p}} &&& \V
\ar @/_1.5pc/[lll]_{\mathsf{l}} \ar
 @/^1.5pc/[lll]_{\mathsf{r}}
 } \ \ \ \ \xymatrix@C=0.3cm{
\V \ar[rrr]^{\mr} &&& \T \ar[rrr]^{\map}  \ar @/_1.5pc/[lll]_{\me}  \ar
 @/^1.5pc/[lll]_{\mr^1} &&& \U
\ar @/_1.5pc/[lll]_{\mathsf{i}} \ar
 @/^1.5pc/[lll]_{\map^1}
 } \ \ \ \
\xymatrix@C=0.3cm{
\U \ar[rrr]^{\map^1} &&& \T \ar[rrr]^{\mr^1}  \ar @/_1.5pc/[lll]_{\map}  \ar
 @/^1.5pc/[lll]_{\mathsf{p}^2} &&& \V
\ar @/_1.5pc/[lll]_{\mr} \ar
 @/^1.5pc/[lll]_{\mathsf{r}^2}
 }
\]
From the bijection between recollements of triangulated categories and \textsf{TTF}-triples, we obtain that $(\ml(\V),\mi(\U),\mr(\V))$, $(\mi(\U),\mr(\V),\map^1(\U))$ and $(\mr(\V),\map^1(\U),\mr^2(\V))$ are \textsf{TTF}-triples in $\T$.

(v) $\Longrightarrow$ (i): From assumption we have the torsion pairs $(\ml(\V),\mi(\U))$, $(\mi(\U), \mr(\V))$, $(\mr(\V),\map^1(\U))$ and $(\map^1(\U),\mr^2(\V))$ in $\T$.
In particular, each of these torsion pairs gives rise to the following diagrams (for the notation see the text before diagram $(\ref{recolarisingfromTTF}))\colon$
\[
(\ml(\V),\mi(\U))\colon \ \xymatrix@C=0.3cm{
\ml(\V) \ar[rrr]^{\mi_{\ml(\V)}} &&& \T  \ar @/^1.5pc/[rrr]^{\mathsf{L}_{\mi(\U)}}  \ar
 @/^1.5pc/[lll]^{\mathsf{R}_{\ml(\V)}} &&& \mi(\U) \ar[lll]_{\mi_{\mi(\U)}} \ \ \ (1)  } \ \ \ \ \ \ \ \ \ \
(\mi(\U),\mr(\V))\colon \ \xymatrix@C=0.3cm{
\mi(\U) \ar[rrr]^{\mi_{\mi(\U)}} &&& \T  \ar @/^1.5pc/[rrr]^{\mathsf{L}_{\mr(\V)}}  \ar
 @/^1.5pc/[lll]^{\mathsf{R}_{\mi(\U)}} &&& \mr(\V) \ar[lll]_{\mi_{\mr(\V)}} \ \ \ (2)  }
\]
\[
(\mr(\V),\map^1(\U))\colon \ \xymatrix@C=0.3cm{
\mr(\V) \ar[rrr]^{\mi_{\mr(\V)}} &&& \T  \ar @/^1.5pc/[rrr]^{\mathsf{L}_{\map^1(\U)}}  \ar @/^1.5pc/[lll]^{\mathsf{R}_{\mr(\V)}} &&& \map^1(\U) \ar[lll]_{\mi_{\map^1(\U)}} \ (3)  } \ \ \
(\map^1(\U),\mr^2(\V))\colon \ \xymatrix@C=0.3cm{
\map^1(\U) \ar[rrr]^{\mi_{\map^1(\U)}} &&& \T  \ar @/^1.5pc/[rrr]^{\mathsf{L}_{\mr^2(\V)}}  \ar @/^1.5pc/[lll]^{\mathsf{R}_{\map^1(\U)}} &&& \mr^2(\V) \ar[lll]_{\mi_{\mr^2(\V)}} \  (4)  }
\]
We now show how from the above diagrams we obtain a ladder $(\U, \T, \V)$  of height three going downwards. From the diagrams $(1)$ and $(2)$ the inclusion functor $\mi_{\mi(\U)}\colon \mi(\U)\lxr \T$ has a left and right adjoint. Then from Lemma~\ref{inducedadjoints} we obtain the recollement $(\mi(\U), \T, \T/\mi(\U))$. In particular, from \cite[Proposition 2.6 (vi), Chapter I]{BR} and Lemma~\ref{inducedadjoints} we obtain the following recollement of triangulated categories$\colon$
\begin{equation}
\label{recollementdiagramone}
\xymatrix@C=0.5cm{
\mi(\U) \ar[rrr]^{\mi_{\mi(\U)}} &&& \T \ar[rrr]^{\mathsf{L}_{\mr(\V)}}  \ar @/_1.5pc/[lll]_{\mathsf{L}_{\mi(\U)}}  \ar
 @/^1.5pc/[lll]_{\mathsf{R}_{\mi(\U)}} &&& \mr(\V)
\ar@{-->} @/_1.5pc/[lll]_{} \ar
 @/^1.5pc/[lll]_{\mi_{\mr(\V)}}
 }
\end{equation}
On the other hand, from diagrams $(2)$ and $(3)$ the inclusion functor $\mi_{\mr(\V)}\colon \mr(\V)\lxr \T$ has a left and right adjoint. Using again \cite[Proposition 2.6, Chapter I]{BR} and Lemma~\ref{inducedadjoints}, but now for the torsion pair $(\mr(\V),\map^1(\U))$, we obtain the following recollement$\colon$
\begin{equation}
\label{recollementdiagratwoandthree}
\xymatrix@C=0.5cm{
\mr(\V) \ar[rrr]^{\mi_{\mr(\V)}} &&& \T \ar[rrr]^{\mathsf{L}_{\map^1(\U)}}  \ar @/_1.5pc/[lll]_{\mathsf{L}_{\mr(\V)}}  \ar
 @/^1.5pc/[lll]_{\mathsf{R}_{\mr(\V)}} &&& \map^1(\U)
\ar@{-->}@/_1.5pc/[lll]_{} \ar
 @/^1.5pc/[lll]_{\mi_{\map^1(\U)}}
 }
\end{equation}
Moreover, from diagram $(4)$ we have the adjoint pair $(\mi_{\map^1(\U)},\mathsf{R}_{\map^1(\U)})$ and Lemma~\ref{inducedadjoints} implies that the functor $\mathsf{R}_{\mr(\V)}$ has a right adjoint, say $\mathsf{R}^1_{\mr(\V)}$. Hence, we have the adjoint triple $(\mathsf{i}_{\mr(\V)}, \mathsf{R}_{\mr(\V)}, \mathsf{R}^1_{\mr(\V)})$. From Lemma~\ref{inducedadjoints} we infer that the recollement $(\ref{recollementdiagramone})$ admits a ladder of height three going downwards.

(iv) $\Longrightarrow$ (vi): Since $(\U, \T, \V)$ restricts to $(\U^{\mathsf{c}}, \T^{\mathsf{c}}, \V^{\mathsf{c}})$, we get the \textsf{TTF}-triple $(\ml(\V^{\mathsf{c}}), \mi(\U^{\mathsf{c}}), \mr(\V^{\mathsf{c}}))$ in $\T^{\mathsf{c}}$. We claim that $\ml(\V^{\mathsf{c}})=\ml(\V)\cap\T^{\mathsf{c}}$, $\mi(\U^{\mathsf{c}})=\mi(\U)\cap\T^{\mathsf{c}}$ and $\mr(\V^{\mathsf{c}})=\mr(\V)\cap\T^{\mathsf{c}}$. We first show that $\mr(\V^{\mathsf{c}})=\mr(\V)\cap\T^{\mathsf{c}}$. Let $\mr(X)$ be an object in $\mr(\V^{\mathsf{c}})$, i.e. $X$ lies in $\V^{\mathsf{c}}$. Since the functor $\mr$ preserves compact objects, the object $\mr(X)$ belongs to $\mr(\V)\cap\T^{\mathsf{c}}$. Conversely, if we take an object $Y$ in $\mr(\V)\cap\T^{\mathsf{c}}$, then $Y=\mr(Y')$ for some $Y'$ in $\V$ and $\mr(Y')$ is compact. Since the counit $\me(\mr(Y'))\lxr Y'$ is an isomorphism and the functor $\me\colon \T\lxr \V$ preserves compact objects, we get that the object $Y'$ lies in $\V^{\mathsf{c}}$. Hence the object $\mr(Y')$ lies in $\mr(\V^{\mathsf{c}})$. Similarly we show that $\ml(\V^{\mathsf{c}})=\ml(\V)\cap\T^{\mathsf{c}}$ and $\mi(\U^{\mathsf{c}})=\mi(\U)\cap\T^{\mathsf{c}}$. We infer that $(\ml(\V)\cap\T^{\mathsf{c}}, \mi(\U)\cap\T^{\mathsf{c}}, \mr(\V)\cap\T^{\mathsf{c}})$ is a \textsf{TTF}-triple in $\T^{\mathsf{c}}$.

(vi) $\Longrightarrow$ (iv): Assuming that  $(\ml(\V)\cap\T^{\mathsf{c}}, \mi(\U)\cap\T^{\mathsf{c}}, \mr(\V)\cap\T^{\mathsf{c}})$ is a \textsf{TTF}-triple in $\T^{\mathsf{c}}$, we show that the recollement
$(\U, \T, \V)$ restricts to a recollement $(\U^{\mathsf{c}}, \T^{\mathsf{c}}, \V^{\mathsf{c}})$, equivalently, we show that $(\ml(\V^{\mathsf{c}}), \mi(\U^{\mathsf{c}}), \mr(\V^{\mathsf{c}}))$ is a \textsf{TTF}-triple in $\T^{\mathsf{c}}$. It suffices to prove that $\ml(\V)\cap\T^{\mathsf{c}}=\ml(\V^{\mathsf{c}})$. Clearly, if $Y$ lies in $\ml(\V^{\mathsf{c}})$, that is, $Y=\ml(V)$ for some $V$ in $\V^{\mathsf{c}}$, then $\ml(V)$ is a compact object in $\T$ since the functor $\ml$ preserves compact objects. Thus, $\ml(\V^{\mathsf{c}})\subseteq \ml(\V)\cap\T^{\mathsf{c}}$. On the other hand, let $X$ be an object in  $\ml(\V)\cap\T^{\mathsf{c}}$. This means that $X$ is of the form $\ml(V)$ for some $V$ in $\V$ and we claim that $V$ lies in $\V^{\mathsf{c}}$.  Indeed, the commutativity of the diagram
\[
\xymatrix{
\Hom_{\V}(V,\coprod_{i\in I}V_i) \ar[d]^{\iso}_{\ml \ \text{f.f.}} \ar[rrrr]^{} & & &  &  \coprod_{i\in I}\Hom_{\V}(V,V_i) \ar[d]^{\iso}_{\ml \ \text{f.f.}}   \\
\Hom_{\T}(\ml(V),\ml(\coprod_{i\in I}V_i)) \ar[rr]^{\iso}_{\ml \ \text{left adjoint}} && \Hom_{\T}(\ml(V),\coprod_{i\in I}\ml(V_i)) \ar[rr]^{\iso}_{\ml(V)\in \T^{\mathsf{c}}} && \coprod_{i\in I}\Hom_{\T}(\ml(V),\ml(V_i))     }
\]
implies that $V$ is compact. Hence, $\ml(\V)\cap\T^{\mathsf{c}}\subseteq \ml(\V^{\mathsf{c}})$ and this completes the proof.
\end{proof}
\end{thm}

We continue with a characterisation of when a recollement of compactly generated triangulated categories admits a ladder of height $n$ going upwards. The proof is analogous to the proof of Theorem~\ref{thmcompactlygenladder} using now Proposition~\ref{propheighttwogoingup} (and its proof), so it is left to the reader.

\begin{thm}
\label{thmcompactlygenladderdual}
Let $(\U, \T, \V)$ be a recollement of compactly generated triangulated categories and $n\geq 3$ a positive odd integer. The following statements are equivalent:
\begin{enumerate}
\item The recollement $(\U, \T, \V)$ admits a ladder of height $n$ going upwards:
\[
\xymatrix@C=0.5cm{
\U  \ar @/^3.0pc/[rrr]^{\mq^1} \ar @/^6.0pc/[rrr]_{}^{\mq^{n-1}} \ar[rrr]^{\mi} &&& \T   \ar @/^3.0pc/[rrr]^{\ml^1} \ar @/_4.5pc/[lll]_{\vdots}  \ar[rrr]^{\me } \ar @/_1.5pc/[lll]_{\mq} \ar @/^6.0pc/[rrr]_{}^{\ml^{n-1}}  \ar @/^1.5pc/[lll]_{\map} &&& \V
\ar @/_1.5pc/[lll]_{\ml} \ar @/_4.5pc/[lll]_{\vdots} \ar
 @/^1.5pc/[lll]_{\mr}   }
 \]

\item There are sequences of functors  $\mq^1, \mq^3, \ldots, \mq^{n-2}\colon \U\lxr \T$ and $\mq,\mq^2, \ldots, \mq^{n-1}\colon \T\lxr \U$ such that $(\mq^{n-1},\mq^{n-2}), \ldots, (\mq^1,\mq)$ are adjoint pairs. In particular, the functors $\mq^1, \mq^3, \ldots, \mq^{n-2}\colon \U\lxr \T$ and $\mq,\mq^2, \ldots, \mq^{n-1}\colon \T\lxr \U$ preserve compact objects.

\item There are sequences of functors $\ml^1, \ml^3, \ldots, \ml^{n-2}\colon \T\lxr \V$ and $\ml ,\ml^2, \ldots, \ml^{n-1}\colon \V\lxr \T$ such that $(\ml^{n-1},\ml^{n-2}), \ldots, (\ml^1,\ml)$ are adjoint pairs. In particular, the functors $\ml^1, \ml^3, \ldots, \ml^{n-2}\colon \T\lxr \V$ and $\ml ,\ml^2, \ldots, \ml^{n-1}\colon \V\lxr \T$ preserve compact objects.

\item The recollement $(\U, \T, \V)$ induces a recollement $(\U^{\mathsf{c}}, \T^{\mathsf{c}}, \V^{\mathsf{c}})$ which admits a ladder of height $n-2$ going upwards.

\item There is sequence of triangulated subcategories $\ml^{n-1}(\V)$, $\mq^{n-2}(\U)$, $\ml^{n-3}(\V)$, $\ldots, \ml^{2}(\V)$, $\mq^{1}(\U)$, $\ml(\V)$, $\mi(\U)$, $\mr(\V)$ in $\T$ such that
\[
\big(\ml^{n-1}(\V), \mq^{n-2}(\U), \ml^{n-3}(\V)\big), \ldots, \big(\ml^2(\V), \mq^1(\U), \ml(\V)\big), \big(\ml(\V),\mi(\U),\mr(\V)\big)
\]
are \textnormal{\textsf{TTF}}-triples in $\T$.

\item There is a sequence of  \textnormal{\textsf{TTF}}-triples in $\T^{\mathsf{c}}$ of the form$\colon$$(\ml^{n-1}(\V)\cap \T^{\mathsf{c}}, \mq^{n-2}(\U)\cap \T^{\mathsf{c}}, \ml^{n-3}(\V)\cap \T^{\mathsf{c}}), \ldots$, $(\mq^{3}(\U)\cap \T^{\mathsf{c}}, \ml^{2}(\V)\cap \T^{\mathsf{c}}, \mq^1(\U)\cap \T^{\mathsf{c}})$, $(\ml^{2}(\V)\cap \T^{\mathsf{c}}, \mq^1(\U)\cap \T^{\mathsf{c}}, \ml(\V)\cap \T^{\mathsf{c}})$.
\end{enumerate}
\end{thm}

As a consequence of Theorem~\ref{thmcompactlygenladder} and Theorem~\ref{thmcompactlygenladderdual} we obtain the following bijections (up to equivalence).

\begin{cor}
\label{corbijections}
For a positive integer $n\geq 1$ we have the following bijections$\colon$
\[
\bigg\{ \stackrel[\displaystyle{\text{of height n going downwards}}]{\displaystyle {(\U, \T, \V) \ \text{is a ladder of compactly generated triangulated categories }}}{\,} \bigg\}
\longleftrightarrow \bigg\{ \X_1, \ldots, \X_{n+2} \stackrel{\Delta \text{\textsf{ed subcat}}}{\subseteq} \T \ |
\]
\[
 \ \ \ \ \ \ \ \ \ \ \ \ \ \ \ \ \ \ \ \ \ \ \ \ \ \ \ \ \ \ \ \ \ \ \ \ \ \ \ \ \ \ \ \ \ \ \ \ \ \  \  \ \ \ \ \ \ \ \ (\X_1,\X_2,\X_3), (\X_2,\X_3,\X_4), \ldots \ (\X_{n},\X_{n+1},\X_{n+2}): \ \textnormal{\textsf{TTF}-triples in} \ \T \bigg\}
\]
\[
\bigg\{ \stackrel[\displaystyle{\text{of height n going upwards}}]{\displaystyle{(\U, \T, \V) \ \text{is a ladder of compactly generated triangulated categories}}}{\,} \bigg\}
\longleftrightarrow \bigg\{ \X_1, \ldots, \X_{n+2} \stackrel{\Delta \text{\textsf{ed subcat}}}{\subseteq} \T \ |
\]
\[
 \ \ \ \ \ \ \ \ \ \ \ \ \ \ \ \ \ \ \ \ \ \ \ \ \ \ \ \ \ \ \ \ \ \ \ \ \ \ \ \ \ \ \ \ \ \ \ \ \ \  \ \ \ \ \ \ \ (\X_{n+2},\X_{n+1},\X_n), (\X_{n+1},\X_n,\X_{n-1}), \ldots \ (\X_{3},\X_{2},\X_{1}): \ \textnormal{\textsf{TTF}-triples in} \ \T \bigg\}
\]
\end{cor}

Note that in \cite{BK} and \cite{IKM:polygons} the authors have also considered sequences of triangulated subcategories such that each two form a stable t-structure.

\begin{rem}
\label{remrelationwithBRttftriples}
Let $\T$ be a compactly generated triangulated category and fix a generating set $\mathcal{S}$ of compact generators in $\T$. Assume that there is a  hereditary torsion pair $(\X,\Y)$ of finite type in $\T$ \cite{BR}, that is, $(\X,\Y)$ is a stable t-structure and $\Y$ is closed under coproducts. Then from \cite[Proposition 1.1, Chapter IV]{BR}  it follows that there is a \textsf{TTF}-triple $(\X,\Y,\Z)$ in $\T$ where $\Z=\Y^{\bot}$. Moreover, the triangulated category $\Y$ is compactly generated by the set $\mathsf{L}_{\Y}(\mathcal{S})$. Hence, we obtain the recollement  $(\Y, \T, \X)$, see diagram $(\ref{recolarisingfromTTF})$.  From \cite[Lemma 1.2, Chapter III]{BR} the functor $\mathsf{R}_{\X}\colon \T\lxr \X$ preserves corpoducts. Then in  \cite[Proposition 1.11, Chapter IV]{BR} the authors provide necessary and sufficient conditions for the hereditary torsion pair $(\X,\Y)$ in $\T$ to induce a torsion pair in $\T^{\mathsf{c}}$. However, these equivalent conditions characterize exactly when the recollement diagram  $(\ref{recolarisingfromTTF})$ admits a ladder of height two going downwards. Moreover, in this case the authors show that the torsion pair $(\X,\Y)$ is compactly generated by the set of objects $\mathsf{R}_{\X}(\mathcal{S})$ and $\X^{\mathsf{c}}=\X\cap \T^{\mathsf{c}}$. This means precisely that $\X={^\bot\Y}$ and $\Y=\big\{\mathsf{R}_{\X}(\mathcal{S})[n]  \ | \  n\geq 0 \big\}$, see \cite[Defintion 2.4, Chapter III]{BR}. Hence, in this way we obtain a description of the sequence of subcategories in Corollary~\ref{corbijections} that provide us the ladder of a recollement of compactly generated triangulated categories.
\end{rem}

\section{Infinite ladders and preprojective algebras}
\label{Section4laddersmono}

In this section we show that the derived category of the preprojective algebra $\Pi_n(\Lambda, Q)$ admits an infinite ladder of period four which restricts to $\mathsf{D}^{\mathsf{b}}(\mathsf{mod})$ and $\mathsf{K}^{\mathsf{b}}(\proj)$. We start by recalling some basic facts about differential graded algebras from \cite{Keller:derivingdgcat, Keller:dgcat}.

\subsection{\bf Differential graded algebras} \ Recall that if $A$ is a dg algebra over a field $k$, then we can construct the homotopy category $\mathsf{H}(A)$ which is a triangulated category and the derived category $\mathsf{D}(A)$ is the localization of $\mathsf{H}(A)$ with respect to the quasi-isomorphisms. Recall from \cite[Theorem 5.3]{Keller:derivingdgcat} that the full subcategory of compact objects of $\mathsf{D}(A)$ coincides with the category $\mathsf{per}(A)$ of perfect $\mathsf{dg}$ $A$-modules. The latter is the smallest full triangulated subcategory
of $\mathsf{D}(A)$ containing $A$ and  closed under finite coproducts and direct summands. We denote by $\mathsf{D}_{\mathsf{fd}}(A)$ the full subcategory of $\mathsf{D}(A)$ consisting of $\mathsf{dg}$ $A$-modules whose total cohomology is finite dimensional, i.e. $\mathsf{D}_{\mathsf{fd}}(A)=\{X \in \mathsf{D}(A) \ | \ \oplus_{n\in \mathbb{Z}}\mathsf{H}^n(X) \ \text{is finite dimensional} \}$.

To proceed, we need the following auxiliary results.

\begin{lem}
\label{fd}
$\mathsf{D}_{\mathsf{fd}}(A)=\{X\in \mathsf{D}(A) \ | \ \oplus_{n\in \mathbb{Z}}\Hom_{\mathsf{D}(A)}(P, X[n]) \ \text{is finite dimensional for any} \ P\in \mathsf{per}(A) \}$.
\begin{proof}
Since we have the isomorphisms $\Hom_{\mathsf{D}(A)}(A,X[i])\iso \Hom_{\mathsf{H}(A)}(A,X[i])\iso \mathsf{H}^0\mathcal{H}om(A,X[i])\iso \mathsf{H}^0(X[i])=\mathsf{H}^i(X)$ and $\mathsf{per}(A)$ is, by definition, closed under direct summands, shifts and extensions, the desired description of $\mathsf{D}_{\mathsf{fd}}(A)$ follows immediately.
\end{proof}
\end{lem}

\begin{lem}
\label{adj}
Let $A$ and $B$ be two $\mathsf{dg}$ algebras. Assume that there is an adjoint pair $(F, G)$ between the derived categories $\mathsf{D}(A)$ and $\mathsf{D}(B)$ such that the functor $F$ restricts to $F\colon \mathsf{per}(A)\lxr \mathsf{per}(B)$. Then the functor $G$ restricts to $G\colon \mathsf{D}_{\mathsf{fd}}(B)\lxr \mathsf{D}_{\mathsf{fd}}(A)$.
\begin{proof}
Let $X$ be an object in $\mathsf{D}_{\mathsf{fd}}(B)$ and let $P$ an object in $\mathsf{per}(A)$. Since there is an isomorphism $\Hom_{\mathsf{D}(A)}(P, G(X)[n])\cong \Hom_{\mathsf{D}(B)}(F(P), X[n])$, for any $n\in \mathbb{Z}$, and $F(P)$ lies in $\mathsf{per}(B)$, it follows from Lemma~\ref{fd} that the hom space $\oplus_{n\in \mathbb{Z}}\Hom_{\mathsf{D}(A)}(P, G(X)[n])$ is finite dimensional. We infer from Lemma~\ref{fd} that the object $G(X)$ lies in $\mathsf{D}_{\mathsf{fd}}(A)$.
\end{proof}
\end{lem}

In the next result we provide a sufficient condition for a recolllement of derived categories of $\mathsf{dg}$ algebras to restrict to $\mathsf{D}_{\mathsf{fd}}$. Compare this result with \cite[Theorem 4.6]{AKLY}.

\begin{prop}
\label{propresticrecotodfd}
Assume that there is a recollement of derived categories of $\mathsf{dg}$ algebras$\colon$
\[
\xymatrix@C=0.5cm{
\mathsf{D}(S)  \ar[rrr]^{\mi} &&& \mathsf{D}(R)  \ar[rrr]^{\me} \ar @/_1.5pc/[lll]_{\mq}   \ar @/^1.5pc/[lll]_{\map } &&& \mathsf{D}(T) \ar @/_1.5pc/[lll]_{\ml} \ar @/^1.5pc/[lll]_{\mr}}  \ \ \ \ \ \ \ \ \ \ \ \ (*)
\]
Consider the following conditions$\colon$
\begin{enumerate}
\item The recollement $(*)$ restricts to a recollement
\[
\xymatrix@C=0.5cm{
\mathsf{D}_{\mathsf{fd}}(S)  \ar[rrr]^{\mi} &&& \mathsf{D}_{\mathsf{fd}}(R)  \ar[rrr]^{\me} \ar @/_1.5pc/[lll]_{\mq}   \ar @/^1.5pc/[lll]_{\map } &&& \mathsf{D}_{\mathsf{fd}}(T) \ar @/_1.5pc/[lll]_{\ml} \ar @/^1.5pc/[lll]_{\mr}}
\]

\item The functor $\ml$ restricts to $\mathsf{D}_{\mathsf{fd}}$ and $\mi(S)$ lies in $\mathsf{per}(R)$.

\item The functor $\mq$ restricts to $\mathsf{D}_{\mathsf{fd}}$ and $\me(R)$ lies in $\mathsf{per}(T)$.

\end{enumerate}
Then \textnormal{(ii)} $\Longrightarrow$ \textnormal{(i)} and \textnormal{(iii)} $\Longrightarrow$ \textnormal{(i)}. If $\map$ preserves coproducts, then all three conditions are equivalent.
\begin{proof}
(ii) $\Longrightarrow$ (i): By Lemma~\ref{adjfromBR} we know that the functors $\mq$ and $\ml$ restrict to $\mathsf{per}(R)$ and $\mathsf{per}(T)$, respectively.
Then by Lemma~\ref{adj} it follows that the functors $\mi$ and $\me$ restrict to $\mathsf{D}_{\mathsf{fd}}$. Since $\mi(S)$ lies in $\mathsf{per}(R)$, i.e. the functor $\mi$ preserves compact objects, Lemma~\ref{adjfromBR} implies that there is an adjoint triple $(\mi,\map,\map^1)$ and therefore an adjoint triple $(\me,\mr,\mr^1)$ by Lemma~\ref{inducedadjoints}. Then the functors $\mi$ and $\me$ restrict to $\mathsf{per}$ and by Lemma~\ref{adj} again we get that the functors $\map$ and $\mr$ restrict to $\mathsf{D}_{\mathsf{fd}}$.
Note that the functor $\ml$ restricts to $\mathsf{D}_{\mathsf{fd}}$ by assumption. It remains to show that $\mq$ restricts to $\mathsf{D}_{\mathsf{fd}}$. Let $X$ be an object in $\mathsf{D}_{\mathsf{fd}}(R)$. From the canonical triangle
\[
\xymatrix{
\ml\me(X) \ar[r]^{} & X \ar[r] & \mi\mq(X) \ar[r] & \ml\me(X)[1] }
\]
and since $\ml\me(X)$ lies in $\mathsf{D}_{\mathsf{fd}}(R)$, we get that $\mi\mq(X)$ lies in $\mathsf{D}_{\mathsf{fd}}(R)$. Since the functor $\mi$ is fully faithful,
there are isomorphisms for all $n\in \mathbb{Z}$ and any $P\in \mathsf{per}(S)\colon$
\[
\Hom_{\mathsf{D}(S)}(P, \mq(X)[n])\cong \Hom_{\mathsf{D}(R)}(\mi(P),\mi\mq(X)[n])
\]
This implies that the space $\oplus_{n\in \mathbb{Z}}\Hom_{\mathsf{D}(S)}(P, \mq(X)[n])$ is finite dimensional and thus Lemma~\ref{fd} shows that $\mq(X)$ lies in $\mathsf{D}_{\mathsf{fd}}(S)$. We infer that $(\mathsf{D}_{\mathsf{fd}}(S), \mathsf{D}_{\mathsf{fd}}(R), \mathsf{D}_{\mathsf{fd}}(T))$ is a recollement.

(i) $\Longrightarrow$ (ii): We only need to check that $\mi(S)$ belongs to $\mathsf{per}(R)$. Since the functor $\map$ preserves coproducts, Lemma~\ref{adjfromBR} implies that $\mi(S)$ lies in $\mathsf{per}(R)$.

The implications (iiii) $\Longrightarrow$ (i) and (i) $\Longrightarrow$ (iii) follow similarly as above.
\end{proof}
\end{prop}

\subsection{Preprojective algebras}\label{subsectionpreproj}
Let $k$ be an algebraically closed field and $Q$ a finite quiver. Denote by
$\overline{Q}$ the double quiver of $Q$ which is obtained from $Q$ by adding for each arrow $a\in Q_1$ an arrow $a^*$ in the opposite direction. Then the preprojective algebra \cite{GelfandPonomarev} is defined as $\Pi(Q):=k\overline{Q}/(c)$ where $k\overline{Q}$ is the path algebra of $\overline{Q}$ over $k$ and $(c)$ is the two-sided ideal generated by $c=\Sigma_{a\in Q_1}(a^{*}a-aa^{*})$. It is known that if $Q$ is Dynkin, then $\Pi(Q)$ is a finite-dimensional selfinjective algebra (see \cite{Reiten}).

{\bf From now on we assume} that the quiver $Q$ is Dynkin of type $\mathbb{A}_n$. For a finite dimensional $k$-algebra $\Lambda$, we denote by $\Pi_n(\Lambda, Q)$ the algebra $\Lambda\otimes_{k}\Pi(Q)$. The latter is a finite dimensional $k$-algebra and is called the path algebra of $\Pi(Q)$ over $\Lambda$, see \cite[subsection 2B]{LZ}. Note that $\Pi_n(\Lambda, Q)$ can be realized as a preprojective algebra of Dynkin species, we refer to \cite{Julian} for more details. The module category of $\Pi_n(\Lambda, Q)$ has objects representations of $\Pi(Q)$ over $\Lambda$, see \cite{ARS}, \cite[Lemma 2.1]{LZ}, \cite[Proposition 4.12]{Julian}. More precisely, for $n=2$ and $n=3$ we have the following descriptions$\colon$
\[
\Mod\Pi_2(\Lambda, Q)=\big\{\xymatrix@C=0.5cm{
X \ar@<-.7ex>[rr]_{f} && Y \ar@<-.7ex>[ll]_-{g}} \ | \ g\circ f=0, \ f\circ g=0 \ \ \text{and} \ X,Y,Z\in \Mod\Lambda \big\}
\]
\[
\Mod\Pi_3(\Lambda, Q)=\big\{\xymatrix@C=0.5cm{
X \ar@<-.7ex>[rr]_{f_1} && Y \ar@<-.7ex>[ll]_-{g_1} \ar@<-.7ex>[rr]_{f_2} && Z \ar@<-.7ex>[ll]_-{g_2} } \ | \ g_1\circ f_1=0=f_2\circ g_2, \ f_1\circ g_1=g_2\circ f_2
\]
\[
\ \ \ \ \ \ \ \ \ \ \ \ \ \ \ \ \ \ \ \ \ \ \ \ \ \ \ \ \ \ \ \ \ \ \ \ \ \ \text{and} \ X,Y,Z\in \Mod\Lambda \big\}
\]
and more generally we have
\[
\Mod\Pi_n(\Lambda, Q)=\big\{\xymatrix@C=0.5cm{
X_{1} \ar@<-.7ex>[rr]_{f_1} && X_{2} \ar@<-.7ex>[ll]_-{g_1} \ar@<-.7ex>[rr]_{f_2} && X_{3} \ar@<-.7ex>[ll]_-{g_2} \ar@<-.7ex>[rr]_{} && \ \cdots \ \ar@<-.7ex>[ll]_-{} \ar@<-.7ex>[rr]_{f_{n-2}} && X_{n-1} \ar@<-.7ex>[ll]_-{g_{n-2}} \ar@<-.7ex>[rr]_{f_{n-1}} && X_{n} \ar@<-.7ex>[ll]_-{g_{n-1}} } \ |
\]
\[
\ \ \ \ \ \ \  \ g_1\circ f_1=0=f_{n-1}\circ g_{n-1}, \ f_{i}\circ g_{i}=g_{i+1}\circ f_{i+1} \, \ \text{for all} \, \ 1\leq i\leq n-2 \ \, \text{and} \ \, X_1, \ldots, X_n\in \Mod\Lambda \big\}
\]
If we restrict to finitely generated $\Lambda$-modules, we get the module category $\smod\Pi_n(\Lambda, Q)$. An object of $\Mod\Pi_n(\Lambda, Q)$ is denoted by $(X_{1},\ldots,X_{n},f_1,g_1,\ldots,$ $ f_{n-1},g_{n-1})$.

We define the following functors$\colon$
\begin{enumerate}
\item The functor $\mt_1\colon \Mod{\Lambda}\lxr \Mod\Pi_n(\Lambda, Q)$ is given by
\[
\xymatrix@C=0.5cm{
X \ar@<-.7ex>[rr]_{1_X} && X \ar@<-.7ex>[ll]_-{0} \ar@<-.7ex>[rr]_{1_X} && X \ar@<-.7ex>[ll]_-{0}\ar@<-.7ex>[rr]_{} && \ \cdots \ \ar@<-.7ex>[ll]_-{} \ar@<-.7ex>[rr]_{1_X} && X \ar@<-.7ex>[ll]_-{0} \ar@<-.7ex>[rr]_{1_X} && X \ar@<-.7ex>[ll]_-{0}}
\]
for a $\Lambda$-module $X$, and given a morphism $a\colon X\lxr X'$ in $\Mod\Lambda$ then $\mt_1(a)=(a,a,\ldots, a)$.

\item The functor $\mt_2\colon \Mod\Lambda\lxr \Mod\Pi_n(\Lambda, Q)$ is given by
\[
\xymatrix@C=0.5cm{
X \ar@<-.7ex>[rr]_{0} && X \ar@<-.7ex>[ll]_-{1_X} \ar@<-.7ex>[rr]_{0} && X \ar@<-.7ex>[ll]_-{1_X}\ar@<-.7ex>[rr]_{} && \ \cdots \ \ar@<-.7ex>[ll]_-{} \ar@<-.7ex>[rr]_{0} && X \ar@<-.7ex>[ll]_-{1_X} \ar@<-.7ex>[rr]_{0} && X \ar@<-.7ex>[ll]_-{1_X} }
\]
for a $\Lambda$-module $X$, and given a morphism $a\colon X\lxr X'$ in $\Mod\Lambda$ then $\mt_2(a)=(a,a,\ldots, a)$.

\item The functor $\mU_{1}\colon \Mod\Pi_n(\Lambda, Q)\lxr \Mod\Lambda$
is given by
\[
\mU_{1}(X_{1},\ldots,X_{n},f_1,g_1,\ldots, f_{n-1},g_{n-1})=X_{1}
\]
on objects, and for a morphism $(a_{1},\ldots, a_{n})$ in $\Mod\Pi_n(\Lambda, Q)$ we have $\mU_{1}(a_{1},\ldots, a_{n})=a_{1}$.

\item The functor $\mU_{2}\colon \Mod\Pi_n(\Lambda, Q)\lxr \Mod\Lambda$
is given by
\[
\mU_{2}(X_{1},\ldots,X_{n},f_1,g_1,\ldots, f_{n-1},g_{n-1})=X_{n}
\]
on objects, and for a morphism $(a_{1},\ldots, a_{n})$ in $\Mod\Pi_n(\Lambda, Q)$ we have $\mU_{2}(a_{1},\ldots, a_{n})=a_{n}$.

\item The functor $\mz_{1}\colon \Mod\Pi_{n-1}(\Lambda, Q)\lxr \Mod\Pi_n(\Lambda, Q)$
is defined by
\[
\mz_1(X_{1},\ldots,X_{n-1},f_{1},g_1,\ldots, f_{n-2},g_{n-2})=(X_{1},\ldots,X_{n-1},0, f_{1},g_1,\ldots,f_{n-2},g_{n-2}, 0, 0)
\]
on objects, and for a morphism $(a_{1}, \ldots, a_{n-1})$ in $\Mod\Pi_2(\Lambda, Q)$ we have $\mz_2(a_{1},\ldots, a_{n-1})=(a_{1},\ldots, a_{n-1}, 0)$.

\item The functor $\mz_{2}\colon \Mod\Pi_{n-1}(\Lambda, Q)\lxr \Mod\Pi_n(\Lambda, Q)$
is defined by
\[
\mz_2(X_{1},\cdots,X_{n-1},f_{1},g_1,\ldots, f_{n-2},g_{n-2})=(0, X_{1},\cdots,X_{n-1},0, 0, f_{1},g_1,\ldots, f_{n-2},g_{n-2})
\]
on objects, and for a morphism $(a_{1}, \ldots, a_{n-1})$ in $\Mod\Pi_{n-1}(\Lambda, Q)$ we have $\mz_2(a_{1},\ldots, a_{n-1})=(0,a_{1},\ldots, a_{n-1})$.
\end{enumerate}

In the next result we show that $\Mod\Pi_n(\Lambda, Q)$ admits a recollement of module categories. The definition of a recollement of module categories is completely analogous to Definition~\ref{defnrecoltriang}. In the abelian case, by definition we have that only the middle functors are exact. For more on recollements of abelian categories see \cite{Pira, Psaroud:homolrecol}.

\begin{prop}
\label{proprecolpreprojalg}
Let $\Lambda$ be a finite dimensional algebra. Then the algebra $\Pi_n(\Lambda, Q)$ admits the following equivalent recollements of module categories$\colon$
\begin{equation}
\label{firstrecollementprepralg}
\xymatrix@C=0.5cm{
\Mod\Pi_{n-1}(\Lambda, Q) \ar[rrr]^{\mz_2} &&& \Mod\Pi_n(\Lambda, Q) \ar[rrr]^{{\mU_1}} \ar @/_1.5pc/[lll]_{}  \ar @/^1.5pc/[lll]_{} &&& \Mod{\Lambda}
\ar @/_1.5pc/[lll]_{{\mt_{1}}} \ar
 @/^1.5pc/[lll]_{\mt_{2}}
 }
\end{equation}
and
\begin{equation}
\label{secondrecollementprepralg}
\xymatrix@C=0.5cm{
\Mod\Pi_{n-1}(\Lambda, Q) \ar[rrr]^{\mz_1} &&& \Mod\Pi_n(\Lambda, Q) \ar[rrr]^{{\mU_2}} \ar @/_1.5pc/[lll]_{}  \ar @/^1.5pc/[lll]_{} &&& \Mod{\Lambda}
\ar @/_1.5pc/[lll]_{{\mt_{2}}} \ar
 @/^1.5pc/[lll]_{\mt_{1}}
 }
\end{equation}
\begin{proof}
We first show that the diagram $(\ref{firstrecollementprepralg})$ is a recollement. From \cite[Remark 2.3]{Psaroud:homolrecol} it suffices to show that $(\mt_1, \mU_1, \mt_2)$ is an adjoint triple with $\mt_1$ (or $\mt_2$) fully faithful and that the kernel $\Ker\mU_1$ is the module category $\Mod\Pi_{n-1}(\Lambda, Q)$. Let $(X'_{1},\ldots, X'_{n},f_1,g_1,\ldots, f_{n-1},g_{n-1})$ be an object in $\Mod\Pi_n(\Lambda, Q)$ and $X$ be an object in $\Mod\Lambda$. Consider a morphism $(a_{1},\ldots,a_{n})\colon \mt_1(X)\lxr (X'_{1},\ldots, X'_{n},f_1,g_1,\ldots, f_{n-1},g_{n-1})$.
Then, it is easy to observe that $a_2=f_1\circ a_1$ and $a_{i+1}=f_i\circ \cdots\circ f_1\circ a_1$ for all $2\leq i\leq n-1$. The assignment $(a_1, f_1\circ a_1, \ldots, f_{n-1}\circ \cdots\circ f_1\circ a_1)\mapsto a_1$
implies that there is a natural isomorphism between $\Hom_{\Pi_n(\Lambda, Q)}(\mt_1(X),(X'_{1},\ldots, X'_{n},f_1,g_1,\ldots, f_{n-1},g_{n-1}))$ and $\Hom_{\Lambda}(X,X_1')$, proving that $(\mt_1,\mU_1)$ is an adjoint pair.
Similarly, if $(a_{1},\ldots,a_{n})\colon (X'_{1},\ldots, X'_{n},f_1,g_1,\ldots, f_{n-1},g_{n-1})\lxr \mt_2(X)$ is a morphism in $\Mod\Pi_n(\Lambda, Q)$, then we get that $a_2=a_1\circ g_1$ and $a_{i+1}=a_1\circ g_1\circ\cdots\circ g_i$ for all $2\leq i\leq n-1$. This shows that the assignment $(a_1, a_1\circ g_1, \ldots, a_1\circ g_1\circ\cdots\circ g_{n-1})\mapsto a_1$ induces a natural isomorphism $\Hom_{\Pi_n(\Lambda, Q)}((X'_{1},\cdots, X'_{n},f_1,g_1,\cdots, f_{n-1},g_{n-1}),\mt_2(X))$ $\cong \Hom_{\Lambda}(X_1',X)$.
Hence, $(\mU_1, \mt_2)$ is an adjoint pair. Clearly, the functor $\mt_1$ is fully faithful and the adjoint triple $(\mt_1,\mU_1,\mt_2)$ implies that $\mt_2$ is also fully faithful. Moreover, the functor $\mz_2$ is fully faithful and the kernel $\Ker\mU_1$ consists of all objects $(X_{1},\ldots,X_{n},f_1,g_1,\ldots, f_{n-1},g_{n-1})$ such that $\mU_1(X_{1},\cdots,X_{n},f_1,g_1,\cdots, f_{n-1},g_{n-1})$ $=0$, that is, $X_1=0$. This implies that the maps $f_1$ and $g_1$ are also zero. Hence, the kernel $\Ker\mU_1$ consists of all objects of the form $(0, X_{2},\ldots,X_{n}, 0, 0, f_2, g_2, \ldots, f_{n-1},g_{n-1})$ and this subcategory is exactly $\mz_2(\Mod\Pi_{n-1}(\Lambda, Q))$. We infer that $(\Mod\Pi_{n-1}(\Lambda, Q), \Mod\Pi_n(\Lambda, Q), \Mod\Lambda)$ is a recollement and similarly we show that $(\ref{secondrecollementprepralg})$ is a recollement as well.

Finally, we show that the recollements $(\ref{firstrecollementprepralg})$  and $(\ref{secondrecollementprepralg})$ are equivalent. Let $(X_{1},\ldots,X_{n},f_1,g_1,\ldots, f_{n-1},$ $g_{n-1})$ be an object in $\Mod\Pi_n(\Lambda, Q)$. We define the endofunctor $\F\colon \Mod\Pi_n(\Lambda, Q)\lxr \Mod\Pi_n(\Lambda, Q)$ by $\F(X_{1},\ldots,X_{n},f_1,g_1,\ldots, f_{n-1}, g_{n-1})=(X_{n},\ldots,X_{1},g_{n-1},f_{n-1},\ldots, g_{1},f_{1})$ on objects, and given a morphism $(a_{1},\ldots, a_{n})$ then $\F(a_{1},\ldots, a_{n})=(a_{n},\ldots, a_{1})$. It follows easily that the functor $\F$ is an equivalence of categories. Consider the following commutative diagram$\colon$
\[
\xymatrix{
    \Mod\Pi_n(\Lambda, Q) \ar[d]_{\F}^{\simeq} \ar[r]^{ \ \ \mU_1} & \Mod\Lambda \ar[d]^{\iden_{\Mod\Lambda}} \\
   \Mod\Pi_n(\Lambda, Q) \ar[r]^{ \ \ \mU_2} & \Mod\Lambda }
\]
It implies that the functors $\mU_{2}\F$ and $\mU_{1}$ are naturally isomorphic. Thus, from \cite[Definition $4.1$, Lemma $4.2$]{PsaroudVitoria:1} we obtain that the recollements $(\ref{firstrecollementprepralg})$ and $(\ref{secondrecollementprepralg})$ of $\Mod\Pi_n(\Lambda, Q)$ are equivalent.
\end{proof}
\end{prop}

We are now ready to state and prove the main result of this paper.

\begin{thm}
\label{doublemorphisminfiniteladder}
Let $\Lambda$ be a finite dimensional algebra over an algebraically closed field $k$ and let $Q$ denote a Dynkin quiver of type $\mathbb{A}_{n}$. For the algebra $\Pi_n(\Lambda, Q)$ (cf. $\S$4.2) the following statements hold$\colon$
\begin{enumerate}
\item There is an infinite ladder $(\mathsf{D}^{}(\Gamma), \mathsf{D}^{}(\Pi_n(\Lambda, Q)), \mathsf{D}^{}(\Lambda))$ of period four, where $\Gamma$ is a $\mathsf{dg}$ algebra such that $\mathsf{H}^0(\Gamma)\iso \Pi_{n-1}(\Lambda, Q)$.

\item There is an infinite ladder $(\mathsf{D}_{\mathsf{fd}}(\Gamma), \mathsf{D}^{\mathsf{b}}(\smod\Pi_n(\Lambda, Q)), \mathsf{D}^{\mathsf{b}}(\smod\Lambda))$ of period four.

\item There is an infinite ladder $(\mathsf{per}(\Gamma), \mathsf{K}^{\mathsf{b}}(\proj\Pi_n(\Lambda, Q)), \mathsf{K}^{\mathsf{b}}(\proj\Lambda))$ of period four.
\end{enumerate}
\begin{proof}
(i) Consider the idempotent element $e_1\in \Pi_n(\Lambda, Q)$ and the right $\Pi_n(\Lambda, Q)$-module $e_1\Pi_n(\Lambda, Q)$. Clearly, the module $e_1\Pi_n(\Lambda, Q)$ is finitely generated projective. Then from \cite[Proposition 5.2]{BazzoniPavarin}, see also \cite[Proposition 2.10]{KalckYang}, we obtain the following recollement of triangulated categories$\colon$
\begin{equation}
\label{recolpreprojalg}
\xymatrix@C=0.5cm{
\mathsf{D}(\Gamma)  \ar[rrr]^{\mi} &&& \mathsf{D}(\Pi_n(\Lambda, Q))  \ar[rrr]^{\me} \ar @/_1.5pc/[lll]_{\mq}   \ar @/^1.5pc/[lll]_{\map } &&& \mathsf{D}(\End(e_1\Pi_n(\Lambda, Q))) \ar @/_1.5pc/[lll]_{\ml} \ar @/^1.5pc/[lll]_{\mr} }
\end{equation}
where $\Gamma$ is a $\mathsf{dg}$ algebra such that $\mathsf{H}^0(\Gamma)\iso \Pi_{n-1}(\Lambda, Q)$ and $\me=e_1\Pi_n(\Lambda)\otimes_{\Pi_3(\Lambda)}^{\mathbb{L}}-$. We now explain the right part of the above recollement.  First, note that the endomorphism algebra $\End(e_1\Pi_n(\Lambda, Q))$ is isomorphic to $\Lambda$. Since $e_1\Pi_n(\Lambda, Q)$ is projective as a right $\Pi_n(\Lambda, Q)$-module, the functor $\me$ is the derived functor of the exact functor $e_1(-)\colon \Mod\Pi_n(\Lambda, Q)\lxr \Mod\Lambda$ which is left multiplication with the idempotent element $e_1$. The latter functor coincides with the functor $\mU_1\colon \Mod\Pi_n(\Lambda, Q)\lxr \Mod\Lambda$, see Proposition~\ref{proprecolpreprojalg}. Then, from Lemma~\ref{exactadjoint} and Proposition~\ref{proprecolpreprojalg} we get the following recollement$\colon$
\begin{equation}
\label{recolpreprojalgfunctors}
\xymatrix@C=0.5cm{
\mathsf{D}(\Gamma)  \ar[rrr]^{\mi} &&& \mathsf{D}(\Pi_n(\Lambda, Q))  \ar[rrr]^{\mU_1} \ar @/_1.5pc/[lll]_{\mq}   \ar @/^1.5pc/[lll]_{\map } &&& \mathsf{D}(\Lambda) \ar @/_1.5pc/[lll]_{\mt_1} \ar @/^1.5pc/[lll]_{\mt_2} }
\end{equation}
The triangle functors $\mt_1$, $\mU_1$ and $\mt_2$ appeared in $(\ref{recolpreprojalgfunctors})$ are the derived functors of the underlying exact functors at the level of module categories. From the recollement $(\ref{secondrecollementprepralg})$, we have the adjoint pair $(\mU_2,\mt_1)$ between $\Mod{\Pi_n(\Lambda, Q)}$ and $\Mod{\Lambda}$. Since both $\mU_2$ and $\mt_1$ are exact functors, Lemma~\ref{exactadjoint} yields an adjoint pair, still denoted by $(\mU_2,\mt_1)$, between the derived categories $\mathsf{D}(\Pi_n(\Lambda, Q))$ and $\mathsf{D}(\Lambda)$. Then, from Lemma~\ref{inducedadjoints} there is a triangle functor $\mq^1\colon \mathsf{D}(\Gamma)\lxr \mathsf{D}(\Pi_n(\Lambda, Q))$ such that $(\mq^1,\mq)$ is an adjoint pair. This implies that the following diagram is a recollement of triangulated categories$\colon$
\begin{equation}
\label{recolpreprojalgqone}
\xymatrix@C=0.5cm{
\mathsf{D}(\Lambda)  \ar[rrr]^{\mt_1} &&& \mathsf{D}(\Pi_n(\Lambda, Q))  \ar[rrr]^{\mq} \ar @/_1.5pc/[lll]_{\mU_2}   \ar @/^1.5pc/[lll]_{\mU_1 } &&& \mathsf{D}(\Gamma) \ar @/_1.5pc/[lll]_{\mq^1} \ar @/^1.5pc/[lll]_{\mi} }
\end{equation}
So far, the above diagram shows that the recollement $(\ref{recolpreprojalgfunctors})$ admits a ladder of height two going upwards. From the recollement of abelian categories $(\ref{secondrecollementprepralg})$, we also have the adjoint pair $(\mt_2, \mU_2)$. Then we get an induced adjoint pair at the level of derived categories still denoted by $(\mt_2, \mU_2)$. Hence, there is a sequence of triangle functors $\mt_1$, $\mU_1$, $\mt_2$, $\mU_2$, $\mt_1$ such that any two consecutive functors form an adjoint pair between $\mathsf{D}(\Pi_n(\Lambda, Q))$ and $\mathsf{D}(\Lambda)$. Then Theorem~\ref{thmcompactlygenladderdual} yields an infinite ladder for $(\ref{recolpreprojalgfunctors})$ of period four going upwards. Recall that the recollement $(\ref{recolpreprojalgfunctors})$ is considered to be
a ladder of height one, so period four means that the fifth recollement that we obtain is the recollement $(\ref{recolpreprojalgfunctors})$. The same method and Theorem~\ref{thmcompactlygenladder} gives an infinite ladder of period four going downwards. We infer that $\mathsf{D}(\Pi_3(\Lambda, Q))$ admits an infinite ladder of period four as follows$\colon$
\begin{equation}
\label{ladderpreprojective}
\xymatrix@C=0.5cm{
\mathsf{D}(\Gamma)  \ar @/^3.0pc/[rrr]^{\mq^1} \ar @/_3.0pc/[rrr]^{\map^1} \ar[rrr]^{\mi } &&& \mathsf{D}(\Pi_n(\Lambda, Q)) \ar @/_4.5pc/[lll]_{\vdots} \ar @/^4.5pc/[lll]^{\vdots} \ar @/^3.0pc/[rrr]^{\mU_2}  \ar[rrr]^{\mU_1} \ar @/_1.5pc/[lll]_{\mq } \ar @/_3.0pc/[rrr]^{\mU_2}  \ar @/^1.5pc/[lll]_{\map} &&& \mathsf{D}(\Lambda) \ar @/_4.5pc/[lll]_{\vdots}
\ar @/_1.5pc/[lll]_{\mt_1} \ar
 @/^1.5pc/[lll]_{\mt_2} \ar @/^4.5pc/[lll]_{\mt_1}^{\vdots}  }
\end{equation}

(ii) $\&$ (iii) Since we have the adjoint triple $(\mU_2,\mt_1,\mU_1)$, Lemma~\ref{adjfromBR} implies that the functor $\mU_2$ preserves compact objects, i.e. it restricts to the category $\mathsf{per}$.
Then from Lemma~\ref{adj} the functor $\mt_1$ restricts to $\mathsf{D}_{\mathsf{fd}}$. Since $(\mi, \map, \map^1)$ is an adjoint triple, it follows from Lemma~\ref{adjfromBR} that $\mi(\Gamma)$ lies in $\mathsf{per}(\Pi_n(\Lambda, Q))$. Then from Proposition~\ref{propresticrecotodfd} the recollement $(\ref{recolpreprojalgfunctors})$ restricts to $(\mathsf{D}_{\mathsf{fd}}(\Gamma), \mathsf{D}_{\mathsf{fd}}(\Pi_n(\Lambda, Q)), \mathsf{D}_{\mathsf{fd}}(\Lambda))$. Since $\Lambda$ and $\Pi_n(\Lambda, Q)$ are finite dimensional algebras (recall that $Q$ is Dynkin of type $\mathbb{A}_n$), so they are considered as $\mathsf{dg}$ algebras concentrated in degree zero, we have triangle equivalences $\mathsf{D}_{\mathsf{fd}}(\Pi_n(\Lambda, Q))\simeq \mathsf{D}^{\mathsf{b}}(\smod\Pi_n(\Lambda))$ and $\mathsf{D}_{\mathsf{fd}}(\Lambda)\simeq \mathsf{D}^{\mathsf{b}}(\smod\Lambda)$. Then we obtain the recollement $(\mathsf{D}_{\mathsf{fd}}(\Gamma), \mathsf{D}^{\mathsf{b}}(\smod\Pi_n(\Lambda, Q)), \mathsf{D}^{\mathsf{b}}(\smod\Lambda))$. As in part (i) we get an infinite ladder of period four for $\mathsf{D}^{\mathsf{b}}(\smod\Pi_n(\Lambda, Q))$ using now the bounded version of Lemma~\ref{exactadjoint}. Finally, since all the involved functors in $(\ref{ladderpreprojective})$ fit into an adjoint triple, Lemma~\ref{adjfromBR} implies that they restrict to compact objects. Then part (iii) follows from Theorem~\ref{thmcompactlygenladder} and Theorem~\ref{thmcompactlygenladderdual}.
\end{proof}
\end{thm}

We now explain that the ladder of Theorem~\ref{doublemorphisminfiniteladder} is in fact a consequence of the Nakayama functor. Compare this with ladders of recollements arising from algebras of finite global dimension treated in \cite{AKLY}.

\begin{rem}
We keep the notation and assumptions as in Theorem~\ref{doublemorphisminfiniteladder}. For simplicity we consider the case $n=2$ and suppose that $\Lambda$ is a selfinjective algebra. We denote by $D\colon \smod{\Lambda}\lxr \smod{\Lambda^{\op}}$ the usual duality, see \cite{ARS}. Consider the Nakayama functor $\nu=D\Pi_2(\Lambda, Q)\otimes_{\Pi_2(\Lambda, Q)}-\colon \smod\Pi_2(\Lambda, Q)\lxr \smod\Pi_2(\Lambda, Q)$ and its right adjoint $\nu^{-1}=\Hom_{\Pi_2(\Lambda, Q)}(D\Pi_2(\Lambda, Q),-)$. We claim that there is an isomorphism  $\mU_2\cong \nu^{-1}\circ \mU_1\circ \nu$ and $(\nu^{-1}\circ \mU_1\circ \nu, T_1)$ is an adjoint pair.

We now show the first claim. Since $\Lambda$ is a selfinjective algebra,  it follows that the algebra $\Pi_2(\Lambda, Q)$ is also selfinjective. Since the functors $\nu$ and $\nu^{-1}$ are exact, it suffices to prove the desired isomorphism for projective modules. Using the decription of $\Pi_2(\Lambda, Q)$ as a Morita ring, see Remark~\ref{remMoritarings} below, we have from \cite[Proposition 3.1]{GP} that the indecomposable projective modules are of the form $\mt_1(P)=(P, P, \iden_{P}, 0)$ or $\mt_2(P)=(P,P, 0, \iden_{P})$, where $P$ is an indecomposable projective $\Lambda$-module. We also refer to \cite[Proposition 2.4]{LZ} for the general case.
Then we compute that $\nu^{-1}\circ \mU_1\circ \nu(P,P, \iden_{P}, 0)=\nu^{-1}(\mU_1(\nu(P), \nu(P), 0, \iden_{\nu(P)}))=\nu^{-1}(\nu(P))\cong P=\mU_2(P,P,\iden_{P}, 0)$ and similarly for the projective $\mt_2(P)$.
For the second claim, it suffices to check the desired adjunction isomorphism for projective modules. Take a projective $\Pi_2(\Lambda, Q)$-module $E=(E_1,E_2, f_{1}, g_{1})$ and a $\Lambda$-module $X$. Since $\Hom_{\Lambda}(\nu^{-1}\circ \mU_1\circ \nu(E), X) = \Hom_{\Lambda}(E_2, X)\cong \Hom_{\Pi_2(\Lambda, Q)}(E, \mt_{1}(X))$ it follows that $(\nu^{-1}\circ \mU_1\circ \nu, \mt_1)$ is an adjoint pair.

In case that $\Lambda$ is a selfinjective algebra (or when we consider just the preprojective algebra $\Pi_n(Q)$), the above considerations shows that from the adjoint pair $(\mt_1, \mU_1)$ between the module categories we can produce the adjoint pair $(\nu^{-1}\circ \mU_1\circ \nu, T_1)$. The main idea of the proof of Theorem~\ref{doublemorphisminfiniteladder} is to lift the underlying exact functors of $(\ref{firstrecollementprepralg})$ and $(\ref{secondrecollementprepralg})$ to derived categories. The key property is that by Proposition~\ref{proprecolpreprojalg} we have the adjoint triples $(\mt_1, \mU_1, \mt_2)$ and $(\mt_2, \mU_2, \mt_1)$ at the level of module categories. Thus when $\Lambda$ is selfinjective, this sequence of adjoints can be interpreted as we explained above via the Nakayama functor.
\end{rem}

In the next remark we explain why we cannot get an analogue of Theorem~\ref{doublemorphisminfiniteladder} for the preprojective algebra $\Pi_n(Q)$ of Dynkin type different from type $\mathbb{A}_n$.

\begin{rem}
Consider the recollement of module categories $(\ref{firstrecollementprepralg})$ and the functor $\mt_1\colon \Mod\Lambda\lxr \Mod\Pi_n(\Lambda, Q)$. We claim that $\mt_1$ is a homological embedding \cite{Psaroud:homolrecol}, that is, there is an isomorphism $\Ext^n_{\Lambda}(X,Y)\cong \Ext^n_{\Pi_n(\Lambda, Q)}(\mt_1(X),\mt_1(Y))$ for all $X,Y\in \Mod\Lambda$ and $n\geq 0$. From the proof of Theorem~\ref{doublemorphisminfiniteladder} we observed that deriving the recollement situation $(\ref{firstrecollementprepralg})$ we get a fully faithful left adjoint $\mt_1\colon \mathsf{D}(\Lambda)\lxr \mathsf{D}(\Pi_n(\Lambda, Q))$. This functor sends complexes concentrated in degree zero to complexes concentrated to degree zero, since the underlying functor $\mt_1$ is exact. This implies that the functor $\mt_1\colon \Mod\Lambda\lxr \Mod\Pi_n(\Lambda, Q)$ is a homological embedding. It was proved by Marks in \cite[Theorem B]{Marks} that if $F\colon \smod{B} \lxr \smod{A}$ is a homological embedding, where $A$ is a preprojective algebra of Dynkin type, then the algebra $A$ needs to be of type $\mathbb{A}_n$ and the algebra $B$ is Morita equivalent to the field $k$. This result together with the method of proving Theorem~\ref{doublemorphisminfiniteladder} shows that we cannot get a similar ladder for the (bounded) derived category of the preprojective algebra $\Pi_n(Q)$ of Dynkin type different from type $\mathbb{A}_n$. If this was the case then the above argument would imply a homological embedding but this contradicts the result of Marks since this can happen only in the  $\mathbb{A}_n$ case.
\end{rem}

\begin{rem}
Let $\Lambda$ be a finite dimensional selfinjective $k$-algebra over a field $k$. From the recollements $(\ref{firstrecollementprepralg})$ and $(\ref{secondrecollementprepralg})$, we have the adjoint triples $(\mU_2,\mt_1, \mU_1)$ and $(\mU_1, \mt_2, \mU_2)$ between $\Mod{\Pi_n(\Lambda, Q)}$ and $\Mod{\Lambda}$. Since  $\mU_2$, $\mt_1$, $\mU_1$ and $\mt_2$ are exact functors, it follows that $\mU_2$, $\mt_1$, $\mU_1$ and $\mt_2$ preserve projective modules. This implies that we get adjoint triples between the stable categories $\underline{\Mod}\Pi_n(\Lambda, Q)$ and $\underline{\Mod}{\Lambda}$, still denoted by $(\mU_2,\mt_1, \mU_1)$ and $(\mU_1, \mt_2, \mU_2)$. Hence, there is an
an infinite sequence of exact functors between the triangulated categories $\underline{\Mod}\Pi_n(\Lambda, Q)$ and $\underline{\Mod}{\Lambda}$ going upwards and downwards such that any two consecutive functors are adjoint pairs$\colon$
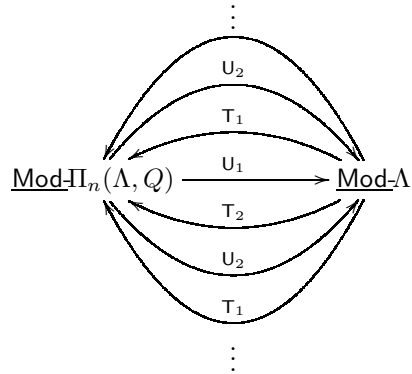
\[
\xymatrix@C=0.5cm{
\underline{\Mod}\Pi_n(\Lambda,Q)  \ar @/^3.0pc/[rrr]^{\mU_2}  \ar[rrr]^{\mU_1}  \ar @/_3.0pc/[rrr]^{\mU_2}  &&& \underline{\Mod}{\Lambda} \ar @/_4.5pc/[lll]_{\vdots}
\ar @/_1.5pc/[lll]_{\mt_1} \ar @/^1.5pc/[lll]_{\mt_2} \ar @/^4.5pc/[lll]_{\mt_1}^{\vdots}  }
\]
Then the kernel $\Ker{\mU_1}$ is a triangulated subcategory of $\underline{\Mod}\Pi_n(\Lambda,Q)$ and as in Theorem~\ref{doublemorphisminfiniteladder} we obtain a periodic infinite ladder
$\mathsf{L}_{\mathsf{tr}}(\Ker\mU_1, \underline{\Mod}\Pi_n(\Lambda,Q), \underline{\Mod}\Lambda)$. Note that $\Ker\mU_1$ is not $\underline{\Mod}\Pi_{n-1}(\Lambda,Q)$ since the functor $\mz_2\colon \Mod\Pi_{n-1}(\Lambda,Q)\lxr \Mod\Pi_n(\Lambda,Q)$ doesn't not preserve projectives.

We mention that since $(\mt_1, \mU_1, \mt_2)$ is an adjoint triple between $\underline{\Mod}\Pi_n(\Lambda, Q)$ and $\underline{\Mod}{\Lambda}$ and $\mU_1$ preserves compact objects, we get immediately from \cite[Theorem 1.7 and 1.9]{Balmer} that there exists an infinite tower of adjoints between $\underline{\Mod}\Pi_n(\Lambda, Q)$ and $\underline{\Mod}{\Lambda}$.
\end{rem}

\begin{exam}
\label{remMoritarings}
Let $\Lambda$ be a finite dimensional algebra and consider the Morita ring (\cite{GP})
\[
\Delta_{(0,0)} = \begin{pmatrix} \Lambda & \Lambda \\ \Lambda & \Lambda \end{pmatrix}
\]
The addition of elements of $\Delta_{(0,0)}$ is componentwise and multiplication is given by
\[
         \begin{pmatrix}
           a & n \\
           m & b \\
         \end{pmatrix}
       \cdot
         \begin{pmatrix}
           a' & n' \\
           m' & b' \\
         \end{pmatrix}=
         \begin{pmatrix}
           aa' & an'+nb' \\
           ma'+bm' & bb' \\
         \end{pmatrix}
\]
The module category of this class of Morita rings was investigated in \cite{GaoPsaroudakis} in connection with aspects of Gorenstein homological algebra. In particular, the module category $\Mod{\Delta_{(0,0)}}$ is equivalent to the double morphism category $\DMor{(\Mod{\Lambda})}$ of $\Mod{\Lambda}$, introduced in \cite[subsection 2.2]{GaoPsaroudakis}. The latter category is exactly the module category $\Mod{\Pi_2(\Lambda, Q)}$ described in the beginning of subsection~\ref{subsectionpreproj}.  Note that the Morita ring $\Delta_{(0,0)}$ is isomorphic to the algebra $\Lambda\otimes_{k} \bigl(\begin{smallmatrix}
k & k \\
k & k
\end{smallmatrix}\bigr)_{(0,0)}$, where the Morita ring $\bigl(\begin{smallmatrix}
k & k \\
k & k
\end{smallmatrix}\bigr)_{(0,0)}$ is the preprojective algebra $\Pi_2(Q)$ over the Dynkin quiver $Q$ of type $\mathbb{A}_2$. Hence, from Theorem~\ref{doublemorphisminfiniteladder} and \cite[Example 2.7]{GaoPsaroudakis} the derived category $\mathsf{D}(\Delta_{(0,0)})$ of the double morphism category admits a periodic infinite ladder $(\mathsf{D}^{}(\Gamma), \mathsf{D}^{}(\Delta_{(0,0)}), \mathsf{D}^{}(\Lambda))$, where $\Gamma$ is a $\mathsf{dg}$ algebra such that $\mathsf{H}^0(\Gamma)\iso \Lambda$, and this ladder restricts to bounded as well as to perfect complexes.
\end{exam}

\end{document}